\documentclass[a4paper,10pt,leqno]{amsart}

\usepackage{amsmath,amsfonts,amssymb}
\usepackage{pifont,graphicx,epsfig}
\usepackage[bf]{caption2}
\usepackage{epstopdf}

\def\C{\hbox{\font\dubl=msbm10 scaled 1000 {\dubl C}}}
\def\Z{\hbox{\font\dubl=msbm10 scaled 1000 {\dubl Z}}}

\def\N{\hbox{\font\dubl=msbm10 scaled 1000 {\dubl N}}}
\def\Q{\hbox{\font\dubl=msbm10 scaled 1000 {\dubl Q}}}
\def\P{\hbox{\font\dubl=msbm10 scaled 1000 {\dubl P}}}
\def\sP{\hbox{\font\dubl=msbm10 scaled 800 {\dubl P}}}

\def\sZ{\hbox{\font\dubl=msbm10 scaled 800 {\dubl Z}}}

\def\Re{{\rm{Re}}\,}
\def\d{{\rm{d}}}

\title[Zeta-Function Theory in the Mathematical Works of Adolf Hurwitz]{Aspects of Zeta-Function Theory in the Mathematical Works of Adolf Hurwitz}

\author[Nicola Oswald, J\"orn Steuding]{Nicola Oswald, J\"orn Steuding}

\date{April 2015}

\begin{document}

\maketitle

\begin{center}
{\it Dedicated to the Memory of Prof. Dr. Wolfgang Schwarz} 
\end{center}
\medskip

\begin{abstract}
Adolf Hurwitz is rather famous for his celebrated contributions to Riemann surfaces, modular forms, diophantine equations and approximation as well as to certain aspects of algebra. His early work on an important generalization of Dirichlet's $L$-series, nowadays called Hurwitz zeta-function, is the only published work settled in the very active field of research around the Riemann zeta-function and its relatives. His mathematical diaries, however, provide another picture, namely a lifelong interest in the development of zeta-function theory. In this note we shall investigate his early work, its origin and its reception, as well as Hurwitz's further studies of the Riemann zeta-function and allied Dirichlet series from his diaries. It turns out that Hurwitz already in 1889 knew about the essential analytic properties of the Epstein zeta-function (including its functional equation) 13 years before Paul Epstein.
\end{abstract}

{\small \noindent {\sc Keywords:} Adolf Hurwitz, Zeta- and $L$-functions, Riemann hypothesis, entire functions}
\bigskip

\tableofcontents

\section{Hildesheim --- the Hurwitz Zeta-Function}

Adolf Hurwitz was born in Hildesheim 26 March 1859 into a Jewish family of merchants. At the age of 17 he had his first publication in mathematics, a joint paper on enumerative geometry with his Hildesheim school teacher Hermann Caesar Hannibal Schubert, a well-known geometer of the time. On promotion of Schubert young Hurwitz started his studies at the polytech university of Munich\footnote{called the {\it K\"oniglich Bayerische Technische Hochschule M\"unchen} which had been renamed to {Technische Universit\"at M\"unchen} in 1970} in 1877 under supervision of the eminent Felix Klein (1849-1925); soon after he moved with Klein to Leipzig where he obtained his doctorate in 1881.
\begin{quote}
''Being a young doctor, in the winterterm 81/82 he moved to the University of Berlin to improve himself by Weierstrass and Kronecker with whom he also had personal contact. Weierstrass (1815-1897) was in particular interested in his doctoral thesis and his function-theoretical efforts, he also proposed the topic for his habilitation thesis.\footnote{''F\"ur das Wintersemester 81/82 bezog er als junger Doktor nochmals die Universit\"at Berlin, um sich bei Weierstra\ss{} und Kronecker, zu denen er auch in pers\"onliche Beziehung trat, noch weiter zu vervollkommnen. Weierstra\ss{} interessierte sich besonders f\"ur seine Dissertation und seine funktionentheoretischen Bestrebungen, auch gab er ihm das Thema f\"ur die Habilitation daselbst...'', \cite{ida}, p. 6; this translation from the original German quotation to English as well as those in the sequel are due to the authors; we have neither corrected the German language to the grammar used nowadays nor have we tried to translate into old-fashioned English. The mentioned university is the {\it Humboldt University}, before 1949 called  {\it Friedrich-Wilhelms-Universit\"at}.}   
\end{quote}
wrote Ida Samuel-Hurwitz, the later wife of Adolf Hurwitz in her dossier \cite{ida}. At Leipzig first obstacles appeared. According to Rowe: ''to pursue an academic career meant applying to habilitate, and the Leipzig regulations required applicants to be graduates of a humanistic gymnasium. Klein must have realized that this stipulation, together with Hurwitz's Jewish confession, virtually ruled out any chance of his habilitating in Leipzig.'', \cite{rowe0}, p. 23.\footnote{Even stronger is a statement of Ida Samuel-Hurwitz about anti-Semitism at Leipzig University during Hurwitz's lifetime: ''Obgleich sein mangelhafter Gesundheitszustand nat\"urlich bekannt war, trachtete man in Deutschland mehrfach, ihn wieder dorthin zu ziehen. Als gleichwertig mit seinem Z\"urcher Lehrstuhl, den nacheinander eine Reihe der hervorragendsten Mathematiker bekleidet hatten, konnten freilich nur Berlin, M\"unchen und Leipzig in Frage kommen. An die letzte Universit\"at war die Berufung eines Juden ausgeschlossen; f\"ur Berlin war er bei Vakanz an erster Stelle vorgeschlagen. Von G\"ottingen war nat\"urlich \"ofters die Rede, doch scheute er einerseits den mathematischen Betrieb dort, w\"ahrend ihm andererseits die Kleinstadt unsympatisch war.'' \cite{ida}, p. 11. English translation: ''Although his defective health was known there were attempts to get him back to Germany. As equivalent to his Zurich chair only Berlin, Munich and Leipzig could be considered. To fill a vacant professorship with a Jew at the latter university was impossible; for Berlin he was recommend for the first place in case of a vacancy. Of course, G\"ottingen was often under discussion, but on one side he disliked the mathematical business there, on the other side he considered the small town unsympathetic.''} Consequently, Hurwitz moved to G\"ottingen for his habilitation in 1882.   


His first thoughts about the Riemann zeta-function and related questions fall into this period. At this time of his career Adolf Hurwitz often visited his nearby hometown Hildesheim. From several letters to his friend from his studies at Munich University (and later renown analyst) Luigi Bianchi (1856-1928) one can get interesting insights about this period. In a letter from 25 May 1882, Adolf Hurwitz had already finished his habilitation and had been working as a private docent in G\"ottingen, he wrote to Bianchi:
\begin{quote}
''Instantly I am visiting home for Pentecost which I can allow myself since Hildesheim can be reached by train from G\"ottingen within 2 hours. (...) In G\"ottingen, I live a pleasant life. 
Prof. H.A. Schwarz is very kind to me and Dr. v. Mangoldt, who is a privat docent for mathematics as well, is a very noble character. Furthermore, in math at G\"ottingen there is Prof. Stern and Prof. Schering.''\footnote{''Augenblicklich bin ich auf Pfingsten zu Besuch zu Hause, was ich mir um so eher erlauben kann, als Hildesheim von G\"ottingen aus in 2 Stunden per Eisenbahn zu erreichen ist. (...) In G\"ottingen f\"uhre ich ein sehr angenehmes Leben. Prof. H.A. Schwarz ist sehr nett gegen mich und Dr. v. Mangoldt, welcher ebenfalls Privatdoc. f\"ur Mathematik ist, ist ein sehr nobler Character. Ferner sind von Mathemat. in G\"ottingen noch Prof. Stern und Prof. Schering.''}
\end{quote}
In particular Hans von Mangoldt (1854-1925) is interesting here. He played an important role in zeta-function theory and its application to questions about prime numbers; his contributions in the 1890s turned out to be essential for the first proof of the prime number theorem in 1896 (details in Section 2). Here, the Riemann zeta-function 
$$
\zeta(s):=\sum_{n\geq 1} n^{-s}=1+2^{-s}+3^{-s}+4^{-s}+\ldots
$$
is building the bridge between primes on the arithmetical side and the powerful tools of complex analysis on the other side. It was young Adolf Hurwitz who invented already in 1881 a simple but important generalization, the nowadays so-called Hurwitz zeta-function, given by
\begin{equation}\label{hurzet}
\zeta(s,\alpha):=\sum_{m\geq 0}(m+\alpha)^{-s},
\end{equation}
where $\alpha$ is a real parameter satisfying $0<\alpha\leq 1$. Actually Hurwitz just considered the case of rational $\alpha$ which is sufficient for applications to number theory; however, for most studies general real numbers $\alpha$ do not cause big trouble. We notice that the series representation from above is absolutely convergent for $\Re s>1$, uniformly in any compact subset. 
\par

We believe that the years 1881-1884 played a very important role in Adolf Hurwitz's life. Of course, he did his doctorate and his habilitation in this period, and in addition he started his important investigations on several topics different to the themes of his doctoral and habilitation theses (on modular forms). With his studies of the Hurwitz zeta-function as well as of real and complex continued fractions he emancipated from his advisor Klein (as well as from his mentor Weierstrass in Berlin) and became an independent mathematician. Concerning continued fractions this is well documented in Oswald \& Steuding \cite{oswald}. In this article we focus on Hurwitz's work on the Riemann zeta-function and its relatives. 

\subsection{Prehistory: From Euler to Riemann}

The infinite series $1+2^{-s}+3^{-s}+\ldots$ defining the zeta-function is linked with the name Leonhard Euler (1707-1783) in many ways. Starting from 1735, he gave various proofs of his famous formula $\zeta(2)={\pi^2\over 6}$; in 1743 he found explicit values for $\zeta(s)$ at the positive even integers \cite{eulerbern}, namely 
$$
\zeta(2k)=(-1)^{k-1}{(2\pi)^{2k}\over 2(2k)!}B_{2k}\qquad\mbox{and}\qquad
\zeta(1-2k)=-{B_{2k}\over 2k}, 
$$
both valid for integral $k$, although Euler had no convergent representation for $\zeta(s)$ when $s$ is a negative real. The numbers $B_m$ are the Bernoulli numbers. Concerning the functional equation linking the values $\zeta(s)$ and $\zeta(1-s)$, Euler claimed that there is no doubt about the truth of this conjecture \cite{eulerfunc}.\footnote{Euler's approach has led to a powerful tool in theoretical physics, called regularization; it was Stephen Hawking \cite{hawk} in 1977 who was the first to use the Riemann zeta-function by this method in order to compute certain quantities in quantum field theory.} Most importantly, although elementary, in 1737 Euler \cite{eulerprod} discovered the analytical analogue of the unique prime factorization in terms of $\zeta(s)$ by observing
$$
\zeta(s)=\sum_{n\geq 1}n^{-s}=\prod_{p\in\sP}(1-p^{-s})^{-1},
$$
valid for $\Re s>1$ (the half-plane of absolute convergence for both, the series and the product), where $\P$ denotes the set of prime numbers. The appearing product is named after Euler and had been used by himself to deduce from the divergent harmonic series (appearing as $\zeta(s)$ with $s\to 1+$) the existence of infinitely many primes by his famous formula
$$
{1\over 2}+{1\over 3}+{1\over 5}+{1\over 7}+\ldots=\log\log \infty
$$ 
(which nowadays is written as $\sum_{p\leq x; p\in\sP}{1\over p}\sim \log\log x$).\footnote{It should be noted that Euler used slightly different notation. For more details on Euler's work we refer to Weil \cite{weil}.} 
\par

The function $\zeta(s)$, however, is named after Bernhard Riemann (1826-1866) who was the first to study the behaviour of $\zeta(s)$ as a function of a complex variable in his path-breaking nine pages article \cite{rie} written at the occasion of his admission to the Prussian Academy of Sciences in 1859. He proved analytic continuation for $\zeta(s)$ to the whole complex plane except for a simple pole at $s=1$ and Euler's conjectured functional equation\footnote{According to Landau \cite{landau}, the onlyother  mathematician to have noticed Euler's work on this topic had been Cahen \cite{cahen2}. This statement is refuted by Malmst\`en \cite{malmsten}, nevertheless some of Euler's works were probably not easily accessible in this time. Only the era of mathematical journals made an easy exchange of mathematical ideas possible.}
$$
\pi^{-{s\over 2}}\Gamma\left({s\over 2}\right)\zeta(s)=\pi^{-{1-s\over 2}}\Gamma\left({1-s\over 2}\right)\zeta(1-s);
$$
Indeed, Riemann gave two rather different proofs of the functional equation, the first one by use of a certain contour integral (which will be the main topic of this section when we shall analyze Hurwitz's paper) and another one based on Poisson's summation formula and Jacobi's transformation formula for the theta-function.\footnote{It was shown by Hamburger \cite{hamburger} that a functional equation of the Riemann-type is equivalent to both, the Poisson summation formula and a theta-function transformation formula.} It follows from the functional equation and basic properties of the gamma-function that the nontrivial zeros (i.e., the non-real zeros) are symmetrically distributed with respect to the real axis and the so-called critical line $\Re s={1\over 2}$. However, Riemann's paper \cite{rie} is probably more important with respect to its not rigorously proved statements and conjectures. Nowadays this very memoir has been considered as a signpost for developing complex analytic methods in order to reveal distribution properties of the prime numbers.


\begin{quote}
''Analytic number theory may be said to begin with the work of Dirichlet, and in particular with Dirichlet's memoir of 1837 on the existence of primes in a given arithmetic progression.''
\end{quote}
This quotation is the first sentence in Davenport's monography on {\it multiplicative number theory} \cite{daven}, refering to the celebrated work of Johannes Peter Gustav Lejeune Dirichlet (1805-1859) on the distribution of primes in arithmetic progressions, resp. prime residue classes around 1837, about one generation before Riemann. In fact, if one seeks for primes in the residue classes modulo $4$, one observes a certain uniform distribution of the primes in the prime residue classes:
$$
\begin{array}{lcccccccl}
4\Z+1\,: & \ldots & 1 & {\bf 5} & 9 & {\bf 13} & {\bf 17} & \ldots &\\
4\Z+2\,: & \ldots & {\bf 2} & 6 & 10 & 14 & 18 & \ldots &\leftarrow\mbox{just one prime}\\
4\Z+3\,: & \ldots & {\bf 3} & {\bf 7} & {\bf 11} & 15 & {\bf 19} & \ldots &\\
4\Z+0\,: & \ldots & 4 & 8 & 12 & 16 & 20 & \ldots &\leftarrow\mbox{no prime at all}
\end{array}
$$
Notice that Euclid's celebrated proof of the infinitude of primes cannot be extended to prove the existence of infinitely many primes in an arbitrary prime residue class $a\bmod\,m$; actually, as proved by Ram Murty (see \cite{murty} for a recent presentation) such an elementary proof is possible if and only if $a^2\equiv 1\bmod\,m$. A first attempt for proving an infinitude of primes in a prime residue class was given by the unfortunate Adrien Marie Legendre \cite{legend}; however, his proof was erroneous, and this was pointed out by Dirichlet \cite{dirichlet37} himself. Before we outline Dirichlet's work on primes in arithmetic progressions we shall briefly mention the circumstances under which Dirichlet wrote his landmark paper \cite{dirichlet37} on this topic (building on Elstrodt's biography \cite{elstrodt}).
\par

His given name {\it Lejeune} means literally {\it the young} and his surname {\it Dirichlet} provides information about the origin of his family being {\it from Richelette}, a small town in Belgium.\footnote{There is something funny about his birthplace: Dirichlet was born in D\"uren, in between Aachen and Cologne, and the house his mother gave birth is {\it Weierstra\ss e 17}; notice that {\it Stra\ss e} is German for {\it street} and the mathematician and contemporary Karl Weierstrass is written {\it Weierstra\ss{}} in German.} His school and university education benefited a lot from the reform of the Prussian educational system implemented by Wilhelm von Humboldt. In 1828, his brother Alexander von Humboldt used his influence to order the young Dirichlet to Berlin.\footnote{see \cite{humbo} for their intensive exchange of letters} Starting from 1831 Dirichlet had a research position at the University of Berlin\footnote{Once again the Humboldt University, formerly Friedrich Wilhems University.} and in the following year he was elected member of the Prussian Academy of Science. In the same year he married Rebecka Mendelssohn, a sister of the composers Fanny Hensel (the grandmother of the number theorist Kurt Hensel) and Felix Mendelssohn Bartholdy.\footnote{The additional surname {\it Bartholdy} had been added to the Jewish surname with respect to the Christian education of the children, a rather common assimilation in these anti-Semitic times.} For the last four years of his short life Dirichlet inherited Carl Friedrich Gauss' chair at the University of G\"ottingen. 
\par

One of his major results is the proof of the existence of infinitely many prime numbers in arbitrary arithmetic progressions $a+m\Z$ with (of course) coprime $a$ and $m$. The main idea is the use of so-called Dirichlet characters, i.e., group homomorphisms $\chi\,:\,(\Z/m\Z)^*\to\C^*$ in order to sift a prime residue class $a\bmod\,m$. Extending these characters to arithmetical functions $n\mapsto \chi(n)$ defined on $\N$ (by letting $\chi(n)=0$ for $n$ not coprime with $m$), 
he defined associated Dirichlet series by
$$
L(s,\chi):=\prod_{p\in\sP} (1-\chi(p)p^{-s})^{-1}=\sum_{n\geq 1}\chi(n)n^{-s}; 
$$
the representations converge for $\Re s>1$ absolutely and the identity between the product and the series is (as in the case of the Riemann zeta-function which appears essentially as $L(s,\chi_0)$ with the principal character $\chi_0$ being constant one) nothing but the unique prime factorization in analytic disguise. Using this in combination with the orthogonality relation for characters, one finds
$$
\sum_{p\equiv a\bmod\,m}p^{-s}\sim {1\over \varphi(m)}\Big(\log\zeta(s)+\sum_{\chi\bmod\,m\atop \chi\neq \chi_0}\overline{\chi}(a)\log L(s,\chi)\Big).
$$
Now the divergence of $\zeta(s)$ for $s\to 1+$ and the regularity and non-vanishing of the values $L(s,\chi)$ with $\chi\neq \chi_0$ at $s=1$ yield the desired statement. 
\par

Dirichlet's approach is rather explicit. In the case of prime moduli $m=p$ his characters are defined by primitive roots and therefore his reasoning for showing $L(1,\chi)\neq 0$ is first limited to such restricted arithmetic progressions. The first appearance of characters can be found in Gauss' work on cyclotomy \cite{gauss}, the modern theory of characters starts with Richard  Dedekind and his supplements to Dirichlet's lectures on number theory \cite{diri}, and the general theory of characters (also for finite nonabelian groups) and group representations is due to Georg Frobenius \cite{frobenius}. For the general and more advanced case of composite moduli Dirichlet built his reasoning on binary quadratic forms $aX^2+bXY+cY^2$. For the number $h(D)$ of equivalence classes of these forms with fixed discriminant $D=b^2-4ac$ he proved the celebrated analytic class number formula \cite{dirichlet39},\footnote{It has been observed by Minkowski \cite{minkowski}, p. 456, that already in Gau\ss' estate one can find a note from 1801 about the class number formula about forty years previous to Dirichlet, however, a proof is missing; moreover, the Disquisitiones \cite{gauss}, p. 369 and p. 466, contain formulae for sums of class numbers without a proof which was first delivered by Siegel \cite{siegel44}.} for example 
$$
h(D)={\omega\over 2\pi}\sqrt{-D}L(1,\chi_D)
$$
for the case $D<0$ of positive definite forms, where $\omega$ denotes the number of units in the ring of integers of $\Q(\sqrt{D})$ and $\chi_D$ is a quadratic character defined by the Kronecker symbol $\chi_D(n)=({D\over n})$. Since $h(D)\geq 1$ this proves $L(1,\chi_D)\neq 0$, and another  subtle reasoning shows that also the non-quadratic characters do not cause any trouble. Dirichlet's reasoning is precise but more complicated than necessary since he used the Jacobi symbol in place of Kronecker's symbol. Anyway, ''Dirichlet's papers, as outstanding for mathematical rigor as Euler had been careless about such matters, ...'' wrote Weil \cite{weil}, p. 7. The first completely analytic proof of the non-vanishing of $L(1,\chi)$ is due to Franz Mertens \cite{mertens} and uses the Dedekind zeta-function associated with quadratic number fields.


Often Schl\"omilch is given credit for the first significant contribution to the analytic theory of $L$-functions beyond Dirichlet's paper. Oscar Xaver Schl\"omilch (1823-1901) was a pupil of Dirichlet; in 1845 he was appointed a professorship at the University of Jena. In 1856, he founded the {\it Zeitschrift f\"ur Mathematik und Physik}, also known as {\it Schl\"omilchs Zeitschrift}, and used this journal to publish several of his results.\footnote{This periodical was only the third one published in Germany after the renowned Crelle journal established in 1826 and run by the Berlin mathematical school, and the less known {\it Archiv f\"ur Mathematik und Physik} founded in 1841 by Johann August Grunert in Greifswald. According to Koch \cite{koch}, the reputation of Schl\"omilch's journal was not very good.} 
\par

In 1849, Schl\"omilch announced in \cite{schloemilch1} the proof of the functional equation 
$$
L(1-s)=\left({2\over \pi}\right)^s\Gamma(s)\sin({\textstyle{\pi s\over 2}})L(s),
$$
where $L(s)=1-3^{-s}+5^{-s}-7^{-s}\pm\ldots$ is the Dirichlet $L$-function to the non-principal character $\chi\bmod\,4$; the proof, however, was published only in 1858 in his journal as \cite{schloemilch2} although the title of his announcement \cite{schloemilch1} is literally 'Exercises for pupils concerning a theorem due to Prof. Dr. Schl\"omilch'. In 1878, he was the first to give a proof of the functional equation for $\zeta(s)$ differing essentially from Riemann's methods \cite{schloemilch278}. Notice that Schl\"omilch was active in educational politics serving from 1874 until 1885 as consultant for the Saxonian Ministery of Education. For more details to his life we refer to Cantor's obituary \cite{cantor}.
\par

However, it was Carl Johan Malmst\`en (1814-1886) being the first to publish a proof of a functional equation for a Dirichlet $L$-function. Interestingly, Malmst\`en's biography shares quite many similarities with Schl\"omilch's. Since 1842 Malmst\`en was professor at the University of Uppsala, and he was among the founders of the renowned journal {\it Acta Mathematica}. He was active as rector of his university as well as privy council and member of the Swedish Reichstag (parliament) for the period 1867-1870. He is famous for his contributions to the Swedish old-age pension insurance and he was member of the Royal Prussian Academy of Sciences as well as the G\"ottingen Academy of Sciences; actually, his daughter Maria Heliodora was married since 1876 to the G\"ottingen mathematician Ernst Christian Julius Schering. In 1849, Malmst\`en \cite{malmsten} published a proof of the functional equation for $L(s)$ independent of Schl\"omilch; moreover he announced that a similar reasoning would prove a corresponding result for 
$$
(1-2^{1-s})\zeta(s)=1-2^{-s}+3^{-s}-4^{-s}\pm\ldots .
$$
Actually, Malmst\`en's paper bears the date 'Upsala 1st May 1846'. It might also be interesting to note that he used already the variable $s$ as Riemann later, which might indicate that Riemann was aware of Malmst\`en's contribution.
\par

But the story before Riemann's memoir appeared is not yet complete. It was Ferdinand Gotthold Max Eisenstein (1823-1852) who gave another proof of a functional equation for a Dirichlet $L$-function and this might have been inspiration for Riemann, too. Eisenstein was another pupil of Dirichlet and member of the Berlin and G\"ottingen Academies of Sciences, the latter one on promotion by Gauss. It was Andre Weil who found in 1989 in Eisenstein's copy of Gauss' {\it Disquisitiones Arithmeticae} on the last and blank page a hand-written proof of the functional equation for $L(s)$ bearing the date 'Scripsi 7 April 1849' (cf. \cite{weil, weil2}), the time when Riemann left Berlin for G\"ottingen. Eisenstein's proof is based on Poisson's summation formula as the second one given by Riemann. It is well-known that Riemann was attending Eisenstein's lectures on elliptic functions in Berlin in the summer 1847 but it seems that their relation was not too good.\footnote{The quotation 'Eisenstein stopped at formal computation' (resp. in German: 'Eisenstein sei bei der formellen Rechnung stehengeblieben') is credited to Riemann; cf. Laugwitz \cite{laugw}, p. 21.} Why did Riemann give two different proofs for the functional equation? Of course, we may only speculate, but a possible answer could be to convince himself about the truth of the result; a further reply could be to stress the importance, and yet another possibility could be awareness of such a result (e.g., Eisenstein's proof).
\smallskip


In 1858 Thomas Clausen (1801-1885) \cite{clausen}, and in 1862 Hermann Kinkelin (1832-1912) \cite{kinkelin} gave alternative proofs for the functional equation for $L(s)$ and generalizations thereof, respectively. Clausen's work \cite{clausen} is related to his studies of Bernoulli numbers. Kinkelin was a Swiss mathematician working essentially in analysis until 1865; later he focused on didactics and mathematics for insurances. In 1865, he was appointed a professorship at the University of Basel\footnote{where around 1896 also Julius Hurwitz, the less famous brother of Adolf Hurwitz, was lecturing for a couple of years; see Oswald \& Steuding \cite{oswald} for details of Julius' interesting life and career.}. Also Kinkelin was politically active (for a conservative party), and in 1886 he was elected Senior Civil Servant (see Schaertlin \cite{schaertlin} for details). His studies \cite{kinkelin} from 1862 include a proof of the functional equation for $L(s,\chi)$ for characters $\chi\bmod\,q$ with $q$ being a prime power or a squarefree odd integer. Kinkelin is not often cited; we can find a note of his work in Landau's {\it Handbuch} \cite{landau09}.\footnote{It is interesting to read Kinkelin's comment on what is nowadays called the Riemann hypothesis: ''The last two equations seem to be satisfied simultaneously only if its root $s$ is of the form $1/2+\tau\sqrt{-1}$; indeed a rigorous proof of this remark has not been given so far (cf. Riemann)''; the German original is: ''Die beiden letzten Gleichungen scheinen nur dann gleichzeitig nebeneinander bestehen zu k\"onnen, wenn ihre Wurzel $s$ von der Form $1/2+\tau\sqrt{-1}$ ist; indessen ist der strenge Beweis dieser Bemerkung noch nicht gelungen (cf. H. Riemann a.a.O.).'' \cite{kinkelin}, p. 18. Though no comment about the arithmetical implications of Riemann's paper are given. It might be that some people of that time considered Riemann's results already as the solution of Gauss' conjecture about the number of primes below a given magnitude.} Kinkelin used the notation $K(s)$ for $\zeta(-s)$ and was very much interested in the relation between Bernoulli numbers and the special values of the zeta-function at the integers. There are further contributions to be mentioned here. Wilhelm Scheibner(1826-1908) \cite{scheibner} was one of the first to analyze Riemann's work; besides he succeeded to establish the case $m=12$. In 1857 Rudolf Lipschitz (1832-1903) \cite{lipschitz} gave yet another proof of the functional equation for $L(s)$, however, we shall consider his results later in another context. 
\smallskip

Riemann \cite{rie} mentioned previous works by Euler and Pafnuty Chebychev (with respect to his elementary approach), but he did not mention Legendre, Eisenstein, Malmst\`en, Schl\"omilch or anybody else. We cannot be sure whether he knew about their results or not. Definitely, Riemann's insight was much deeper and his methods were far beyond what his contemporaries had at hand.

\subsection{The Hurwitz Identity}

Adolf Hurwitz's article \cite{hurwzeta} from 1881 starts as follows:
\begin{quote}
''(...) I achieved general formulae of the same simplicity as those found by Schl\"omilch and Riemann, which appear as particular cases. These general formulae might be of special interest with respect to their relation to functions playing a substantial role in Dirichlet's deep investigations concerning the number of classes of binary quadratic forms and prime numbers in an arithmetic progression.''\footnote{''(...) gelangte ich zu allgemeinen Formeln, die von gleicher Einfachheit sind, wie die von Schl\"omilch und Riemann gegebenen, und dieselben \"uberdies als spezielle F\"alle enthalten. Diese allgemeinen Formeln d\"urften um so mehr ein Interesse beanspruchen, als sie sich auf eben die Functionen beziehen, die bei den tiefen Untersuchungen Dirichlet's \"uber die Anzahl der Klassen bin\"arer quadratischer Formen und \"uber die in einer arithmetischen Reihe enthaltenen Primzahlen eine wichtige Rolle spielen.'', \cite{hurwzeta}, p. 73/74.} 
\end{quote}
Notice that \cite{hurwzeta} was published in {\it Schl\"omilchs Zeitschrift}. Therefore, it is natural to mention Schl\"omilch's contribution besides Riemann and Dirichlet, however, in the whole paper neither references to Malmst\'en, Eisenstein nor Kinkelin can be found.
\par
 
The 22 years old young doctor Adolf Hurwitz signs his paper with the date ''Hildesheim, den 10. Oktober 1881'' although his scientific address at that time was Berlin. The letter exchange from this period with his friend Luigi Bianchi is very interesting for several reasons. They give the picture of two extraordinarily gifted mathematicians who kept their friendship despite the geographical distance between each other.     
\begin{quote}
''Infinitely I have enjoyed to have been welcomed by your letter; your last letter with your photography (the people are stupid, they go and let them photographed) I have received too.''\footnote{''Unendlich habe ich mich gefreut mit einem Briefe von Dir empfangen zu werden; auch Deinen letzten Brief mit Deiner Photographie (Die Leute sind dumm, sie gehen und lassen sich photografiren) habe ich erhalten'', \cite{bianchi}, p. 77.} 
\end{quote}
wrote Hurwitz to Bianchi in a letter with date 23 December 1881, \cite{bianchi}, p. 77; a later letter of Hurwitz to Bianchi from 6 April 1883 must have contained a photograph of Hurwitz, too; \cite{bianchi}, p. 91. They continued writing letters to each other their whole lives, however, it seems they never met again in person. Of particular interest are mathematical topics in their correspondence. 
In a letter to Hurwitz with date 28 August 1881 Bianchi wrote (see Figure 1): 
\begin{quote}
''I am grateful for your function-theoretical note. If I am not wrong the path to generalize these results is in considering Dirichlet series 
$$
\sum \left({D\over n}\right){1\over n^s}
$$ 
where $n$ runs thorugh all prime numbers with respect to $D$. For $D=2$ one gets Schl\"omil(ch)'s series and for $D=3$ one gets yours.''\footnote{''Ich danke Dir f\"ur Deine functionentheoretische Mittheilung. Wenn ich nicht irre liegt der Weg zur Verallgemeinerung dieser Resultate in der Betrachtung der Dirichlet'schen Reihen 
$$
\sum \left({D\over n}\right){1\over n^s}
$$
wo $n$ alle Primzahlen zu $D$ durchlaeuft. F\"ur $D=2$ bekommt man die Schl\"omil(ch)'sche und f\"ur $D=3$ (die Deinige).''}
\end{quote}
Unfortunately, the previous letter of Adolf Hurwitz got lost; however, in \cite{bianchi}, Opere, Vol. XI: Corrispondenza, one can read Hurwitz's reply to Bianchi from 24 September 1881:
\begin{quote}
''Your serendipity has found the correct generalization; however, I have had the same idea a few days after sending you my previous letter on my own, since the easy treatable case of a prime number had necessarily guided me. Today I shall inform you about the complete results...''\footnote{''Dein Sp\"ursinn hat die richtige Verallgemeinerung getroffen; ich hatte allerdings auch ein paar Tage nach der Absendung meines vorigen Briefes dieselbe Idee, da mich der leicht zu behandelnde Fall einer Primzahl mit Nothwendigkeit darauf f\"uhrte. Ich theile Dir nun heutige die fertigen Resultate mit...'', \cite{bianchi}, p. 74.} 
\end{quote}
Luigi Bianchi was a student of Enrico Betti and Ulisse Dini; he received his doctorate already in 1877. Afterwards he continued his studies at the University of Munich where he became a friend of Adolf Hurwitz. Their correspondence gives indeed a lively picture of a deep and lifelong friendship and it seems there was no competition between the two young stars.
\smallskip

We shall briefly outline Hurwitz's proof (using more or less his notation). First of all, given integers $1\leq a\leq m$, he defines  
$$
f(s,a):=a^{-s}+(a+m)^{-s}+(a+2m)^{-s}+\ldots=\sum_{n\equiv a\bmod\,m}n^{-s};
$$
in the now common notation (\ref{hurzet}) these Dirichlet series equals $m^{-s}\zeta(s,a/m)$, and we recover the Riemann zeta-function as $\zeta(s)=f(s,1)$ in the case $m=1$. Building on Riemann's first proof of the functional equation, Hurwitz starts with Euler's integral representation of the gamma-function, 
$$
\Gamma(s)=\int_0^\infty t^{s-1}\exp(-t)\d t,
$$
valid for $\Re s>0$. Substituting $t=nx$ and summing up over $n\equiv a\bmod\,m$, he gets
$$
\Gamma(s)\sum_{n\equiv a\bmod\,m}n^{-s}=\sum_{n\equiv a\bmod\,m}\int_0^\infty x^{s-1}\exp(-nx)\d x.
$$
By uniform convergence, interchanging summation and integration is allowed; the resulting sum of exponentials is a geometric series with value
\begin{equation}\label{thus}
\sum_{n\equiv a\bmod\,m}\exp(-nx)=\sum_{k\geq 0}\exp(-(km+a)x)={\exp(-ax)\over 1-\exp(-mx)}.
\end{equation}
It might be interesting to notice that already in Dirichlet's approach \cite{dirichlet37} one can find a comparable integral representation for the product of a Dirichlet series with the gamma-function. Actually, Hjalmar Mellin \cite{mellin1,mellin2} developed a theory for such integral representations coining the notion of a Mellin transform. It might be interesting to notice that in 1881/82 Mellin studied at the University of Berlin (cf. Lindel\"of \cite{lindeloef}) so he might have met Hurwitz in person when he was starting his investigations of the Hurwitz zeta-function. Previous to Mellin it was Eug\`ene Cahen \cite{cahen2} who studied for his doctorate such transformations, however, his dissertation had been criticised for its many errors by Adolf Hurwitz in his review \cite{hurwreport1} for the Jahresberichte; it was later corrected by Oskar Perron \cite{perron} (who also removed some mistake of Kronecker \cite{kronecker}).

We continue with Hurwitz's proof. it follows from (\ref{thus}) that    
$$
\Gamma(s)f(s,a)=\int_0^\infty {\exp((m-a)x)\over \exp(mx)-1}x^{s-1}\d x.
$$
Now define the integral\footnote{Such integrals are now called Hankel integrals in honour of Riemann's pupil Hermann Hankel who used them in his studies \cite{hankel} of the Gamma-function (although already Riemann \cite{rie} used them earlier).}
$$
J(s)=\int_C\exp((s-1)\log(-x)){\exp((m-a)x)\over \exp(mx)-1}\d x,
$$
where the contour $C$ is the path starting from $+\infty$ along the real axis in the upper half-plane to a real point $A>0$, then counterclockwise continuing on the circle $\circ$ of radius $A$ around the origin and then back to $+\infty$ along the real axis in the lower half-plane; we may assume that $A$ is so small that the only singularity of the integrand inside $C$ is at the origin, and that there is no singularity on $C$. For positive $x$ the imaginary part of the logarithm $\log(-x)$ is equal to $-i\pi$ on the positive part of the real axis and equal to $+i\pi$ on the negative part. Consequently, 
\begin{eqnarray*}
J(s)&=&\int_A^\infty (\exp(-i\pi s)-\exp(i\pi s))x^{s-1}{\exp((m-a)x)\over \exp(mx)-1}\d x+\\
&&+\int_\circ \exp((s-1)\log(-x)){\exp((m-a)x)\over \exp(mx)-1}\d x.
\end{eqnarray*}
The integral over the small circle $\circ$ tends for $\Re s>1$ to zero as $A\to 0$, giving
$$
\lim_{A\to 0}J(s)=-2i\sin(\pi s)\Gamma(s)f(s,a).
$$
Using the formula 
$$
\Gamma(s)\Gamma(1-s)={\pi\over \sin\pi s}, 
$$
Hurwitz deduces 
$$
f(s,a)={i\over 2\pi} \Gamma(1-s)\lim_{A\to 0}J(s).
$$
Since the integral $\lim_{A\to 0}J(s)$ converges uniformly on compact sets it defines an entire function. Consequently, the only possible singularities of $f(s,a)$ are the poles of $\Gamma(1-s)$ at the positive integers; however, since $f(s,a)$ is obviously regular for $s=2,3,\ldots$, it only remains to consider the local behaviour at $s=1$. In view of $\Gamma(z+1)=z\Gamma(z)$ (which was already known to Euler) one has 
$$
\Gamma(1-s)={-1\over s-1}+\ldots
$$
and in combination with
$$
\lim_{A\to 0}J(1)=\int_0^\infty (-x)^{s-1}{\exp((m-a)x)\over \exp(mx)-1}\d x ={2\pi i\over m}
$$
it follows that {\it $f(s,a)$ admits an analytic continuation to $\C$ except a simple pole at $s=1$ with residue $1/m$.} In a similar way Hurwitz deduces that the values $f(s,a)$ are rational functions in $a$ and $m$, e.g.,
$$
f(-2,a)=-{a(m-a)(m-2a)\over 6m}.
$$
In modern terms this formula is a special case of the formula
$$
\zeta(-n,{\textstyle{a\over m}})=-{1\over n+1}B_{n+1}({\textstyle{a\over m}})\quad (\ =m^{-n}f(-n,a)\ ) 
$$
with the Bernoulli polynomials $B_k(X)$ defined by the expansion 
$$
{t\exp(tX)\over \exp(t)-1}=\sum_{k\geq 0}B_k(X){t^n\over n!},
$$
giving $B_3(X)=X^3-{3\over 2}X^2+{1\over 2}X$ in particular. This as well as much of Hurwitz's investigation of the functional equation in modern form can be found in Lang's books \cite{lang0,lang}.
\par

Hurwitz's next aim is the proof of a functional equation for his function $f(s,a)$. For this purpose the integral $J$ is evaluated in a second different way. In place of $C$ where the integral was taken over a circle of radius $A$ centered at the origin the new path is a counterclockwise oriented rectangle with vertices $(2N+1){\pi\over m}(\pm 1\pm i)$ which includes the simple poles of the integrand at $x={2\pi in\over m}$ with integers $n$ such $\vert n\vert\leq N$ and excluding all others. Thus the new integral differs by the old one by the sum over all residues resulting from these poles except the one at the origin. It is no difficulty to compute the residues for $\pm n$ as 
$$
{1\over m}\left({2\pi n\over m}\right)^{s-1}\exp(\mp{\textstyle{i\pi\over 2}}(s-1)\mp 2\pi i{\textstyle{an\over m}}).
$$
Now let $\Re s<0$. Since then the integrals over the edges of the rectangular paths vanish as $N\to\infty$, Hurwitz arrives at the equation
$$
f(s,a)={1\over \pi}\left({2\pi\over m}\right)^s\Gamma(1-s)\sum_{b=1}^m\sum_{k=0}^\infty (mk+b)^{s-1}\cos({\textstyle{\pi\over 2}}(s-1)+2\pi{\textstyle{ab\over m}}).
$$
Finally, after replacing $s$ by $1-s$, this leads to the celebrated {\it Hurwitz formula}
\begin{equation}\label{hurwitzidentity}
f(1-s,a)={1\over \pi}\left({2\pi\over m}\right)^{1-s}\Gamma(s)\sum_{b=1}^m\cos(2\pi{\textstyle{ab\over m}}-{\textstyle{\pi s\over 2}})f(s,b).
\end{equation}
By analytic continuation the latter formula holds for all values of $s$. 
\par

Using this powerful identity Hurwitz investigated next his main target, that is the Dirichlet series
$$
F(s,D)=\sum_{n\geq 1} \left({D\over n}\right)n^{-s},
$$
where $n\mapsto ({D\over n})$ is the Jacobi symbol associated with a squarefree positive integer $D\equiv 1\bmod\,4$ and the summation is over all odd $n$ coprime with $D$; for the cases $D\equiv 2$ or $3\bmod\,4$ a harmless factor has to be multiplied, however, for the sake of brevity we shall consider only $D\equiv 1\bmod\,4$ here. Using a decomposition of $F(s,D)$ into his functions $f(s,a)$ with $m=\vert D\vert$ and exploiting quadratic reciprocity, Hurwitz proves that {\it with the only exception $D=1$ the function $F(s,D)$ is an entire function satisfying the functional equation}
$$
F(1-s,D)=\left({D\over \pi}\right)^{s-{1\over 2}}{\Gamma({\textstyle{s\over 2}})\over \Gamma({\textstyle{{1-s\over 2}}})}F(s,D).
$$
In principle, his reasoning is a simple application of the representation
$$
L(s,\chi\bmod\,m)=\sum_{a\bmod\,m}\chi(a)f(s,a)
$$
in combination with his identity (\ref{hurwitzidentity}) and explicit formulas for Gauss sums to primitive characters $\chi\bmod\,D$, namely 
$$
\sum_{a\bmod\,D}\chi(a)\exp\left({2\pi ia\over D}\right)=i^{{1\over 2}(1-\chi(-1))}D^{1\over 2},
$$
which appear first by Dirichlet \cite{dirichlet37} (continuing Gauss' work \cite{gauss}). Hurwitz's reasoning is actually slightly more complicated as it is nowadays presented in textbooks since the notion of primitive characters was at his time not well formulated. As a matter of fact primitive characters have been used by Edmund Landau in his 1909 monography \cite{landau09}; for some reason (see 'Quellenangaben' in \cite{landau09}, p. 898) his reference for that is an article by Lipschitz \cite{lip89} from 1889 (more precisely he refers to pages 142-144), however, as already pointed out by Narkiewicz \cite{narkie}, the characters considered by Lipschitz are primitive but Lipschitz did not invent the distinction in between primitive and imprimitive characters. Whereas Hadamard did not consider the case of arithmetic progressions in his proof of the prime number theorem de la Vall\'ee Poussin worked with characters, however, there is no notion of a primitive character in his works (since his approach does not rely on the functional equation). A Dirichlet $L$-function to an imprimitive character does not satisfy a functional equation of the Riemann type.
\smallskip

We conclude this paragraph with another quotation from the letter Hurwitz sent to Bianchi on 24 September 1881:
\begin{quote}
''My social life is very active thanks to music. I accompany my singing (female) friends and play with them (you don't need to think anything bad). Recently, we played eight-handed with two pianos: three women and myself. -- I don't know the article of Lipschitz and cannot receive it from here. From Berlin I can tell you about its contents if it is still of any value for you.''\footnote{''Mein gesellschaftliches Leben ist durch die Musik ein sehr reges. Ich begleite meine Freundinnen zum Gesang und spiele mit Ihnen (Du brauchst nichts Schlechtes zu denken). Neulich haben wir sogar 8-h\"andig auf 2 Claviren musicirt: 3 Damen und ich. -- Den Aufsatz von Lipschitz kenne ich nicht und leider kann ich mir denselben hier auch nicht verschaffen. Von Berlin aus kann ich Dir \"uber denselben berichten, wenn es dann f\"ur Dich noch Werth hat.'' \cite{bianchi}, p. 76.} 
\end{quote}
Their correspondence \cite{bianchi} contains two further notes on this topic, the first one is a variation of Dirichlet's class number formula Hurwitz found (and reported to Bianchi in a letter from 23 December 1881) and the second is about the delay in the publishing process (in another letter with date 20 March 1882). Not investigating what Lipschitz had done is juvenile ommission. In the following we shall consider the contributions made by Lipschitz and some other contemporaries. 

\subsection{Aftermath and Further Generalizations by Lipschitz, Lerch, and Epstein}

Rudolf Lipschitz was a student of Dirichlet. In 1864 he became professor in Bonn. He had been considered as one of the leading analysts of his time. In his first note on Dirichlet series \cite{lipschitz} (and one of the first of all his papers) from 1857 Lipschitz studied the real and imaginary part of the $\zeta(s)$ defining Dirichlet series and generalizations thereof, namely 
$$
\sum_{n=-\infty}^{+\infty} \exp(2\pi \sqrt{-1}v)((n+u)\sqrt{-1}+k)^{-\sigma};
$$
he obtained integral representations and identities for these infinite series. Here $\sigma$ is assumed to be a real variable; in that sense Lipschitz gave birth to the funny mixture of writing $s=\sigma+it$ with greek and latin characters for the complex variable (which is often attributed to Landau; the imaginary part $t$ stems from Riemann's paper), only the use of the letter 'i' for the imaginary unit $\sqrt{-1}$ cannot be found in his papers. With his second article \cite{lip89} from 1889 Lipschitz succeeded to give a more or less complete account of Dirichlet $L$-functions within an analysis of a more general Dirichlet series. Following Riemann's first proof, Lipschitz \cite{lip89} obtained the identity
$$
\Gamma\left({1-s\over 2}\right)\sum_{n=-\infty}^{+\infty}{\exp(-2\pi iv(n+u))\over ((n+u)^2\pi)^{s\over 2}}=
\Gamma\left({s\over 2}\right)\sum_{n=-\infty}^{+\infty}{\exp(2\pi inu)\over ((n+v)^2\pi)^{1-s\over 2}},
$$
valid for $\Re s\in(0,1)$. Using Gauss sums, he deduced the functional equation for quite a few Dirichlet $L$-functions; his reasoning is rather general. Later, Vall\'ee Poussin \cite{vapo} found another simplified proof for the functional equation of Dirichlet $L$-functions; the method of proof relies on Riemann's second proof of the functional equation for $\zeta(s)$. For this purpose Vall\'ee Poussin proved a character version of the Jacobi theta-function transformation formula. His proof became standard, appearing first in Landau's {\it Handbuch} \cite{landau09} (see also Davenport \cite{daven}).
\par

Besides Hurwitz, Lipschitz gave credit to Mathias Lerch\footnote{or {\it Maty\'a\v s Lerch} in Czech; Lerch started school only at the age of nine years because of a severe injury of his left leg, a handicap he suffered for the rest of his life. For details of his unfortunate career and difficult character we refer to Porubsk\'y \cite{porubsky}.} (1860-1922) \cite{lerch87} who studied two years earlier the infinite series which is nowadays well-known as the Lerch zeta-function, namely (in modern standard notation)
$$
L(\lambda,\alpha,s)=\sum_{m\geq 0} \exp(2\pi i\lambda m)(m+\alpha)^{-s} 
$$
with real $\lambda$ and $\alpha$ (so that $L(1,\alpha,s)=\zeta(s,\alpha)$ with $\alpha\in(0,1])$. For this function Lerch obtained analytic continuation to the whole complex plane except a simple pole at $s=1$ in the case of integral $\lambda$ and he derived a functional identity of the form
\begin{eqnarray*}
L(\lambda,\alpha,1-s)&=&(2\pi)^{-s}\Gamma(s)\large(\exp({\textstyle{\pi is\over 2}}-2\pi i\alpha\lambda)L(-\alpha,\lambda,s)\\
&&\qquad +\exp(-{\textstyle{\pi is\over 2}}+2\pi i\alpha(1-\lambda))L(\alpha,1-\lambda,s)\large), 
\end{eqnarray*}
valid for $\lambda\in(0,1), \alpha\in(0,1]$, and $\Re s\in(0,1)$. This formula includes the Hurwitz formula as limiting case $\lambda\to 0+$. It is interesting that Lerch even allowed complex values for $\lambda$ and $s$. It seems that complex numbers and methods were widely accepted in Lerch's generation but less in the former. His reasoning is also based on Riemann's first approach to the functional equation. As particular case one deduces a result due to Lipschitz \cite{lipschitz} and Kronecker \cite{kronecker}, namely the identity
$$
\sum_{m=-\infty}^{+\infty}{\exp(2\pi im\lambda)\over m-\alpha}=2\pi i{\exp(2\pi i\alpha\lambda)\over 1-\exp(\pi i\alpha\lambda)}.
$$
Lerch mentioned previous works by Malmst\`en \cite{malmsten} and Lipschitz \cite{lipschitz}; moreover, in a later paper \cite{lerch92} he even called the latter series the {\it Lipschitz function}.\footnote{In 1889, Alfred Jonquiere \cite{jonquiere} obtained the analogue of Hurwitz's formula for the polylogarithm, however, Lerch's result is more general.} 
Malmst\`en investigated the even more general series
$$
\sum_{n=-\infty}^{+\infty}{\exp(2\pi inu)\over ((n+u)^2+v^2)^{s\over 2}},
$$
which has some similarity with the so-called Epstein zeta-function associated with a quadratic form. As a matter of fact, already in a letter to Bianchi from 1881 Hurwitz had mentioned these series: 
\begin{quote}
''In order to gain information about class numbers from these formulae one would have to investigate the corresponding relationship for 
$$
\sum \left({1\over ax^2+2bxy+cy^2}\right)^s
$$
''\footnote{''Um aus diesen Formeln Vortheil f\"ur die Classenzahlen zu gewinnen, m\"usste man den entspr. Zusammenhang f\"ur die 
$$
\sum \left({1\over ax^2+2bxy+cy^2}\right)^s
$$
erforschen.''} 
\end{quote}
This series appeared already in Dirichlet's paper \cite{dirichlet39}. A thorough analysis of such Dirichlet series associated with quadratic forms was realized by Paul Epstein (1871-1939) in his 1903 paper \cite{epstein}.
\par

\begin{figure}[h]
\includegraphics[height=7.7cm]{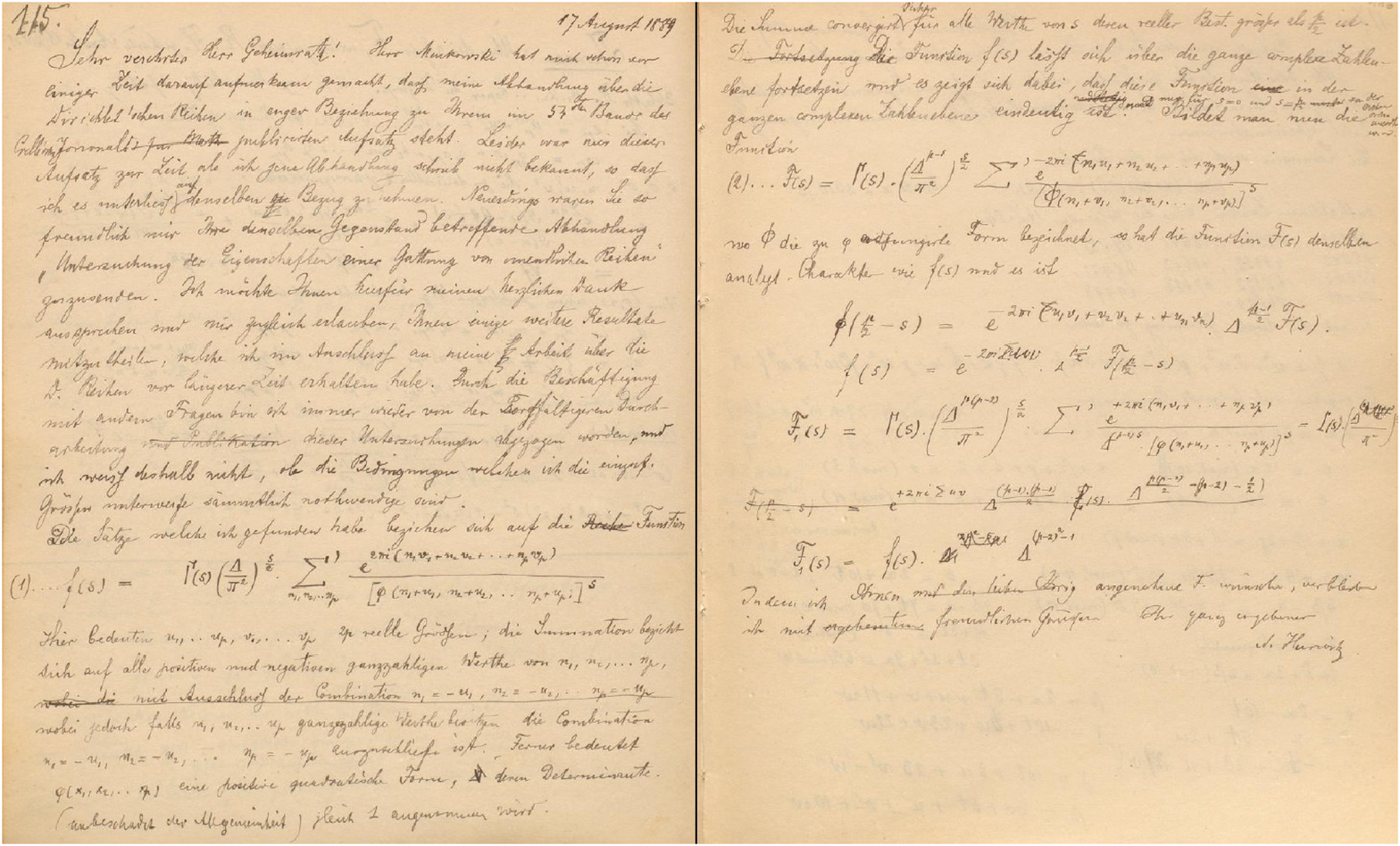}
\caption{The proposal for a letter to Lipschitz with date 17 August 1889, found in Hurwitz's Mathematical Diary no. 6, p. 115/116.}
\end{figure}
In Hurwitz's Diary no. 6 with date 17 August 1889 we can find the following lines which seem to be a proposal for a letter he had probably written to Lipschitz, however, this letter seems to be lost:\footnote{At least, Hurwitz is not mentioned in Scharlau's list \cite{scharlau} of letters Lipschitz had received from his contemporaries. There appears also no correspondence between Lipschitz and Hurwitz in the latter one's estate, neither in Zurich nor in the archive of the Staats- und Universit\"atsbiliothek G\"ottingen where most f his letters are stored.}
\begin{quote}
''Dear privy councillor! Mr Minkowski had informed me some time ago that my treatise on Dirichlet series is intimately related with your article published in the 54th issue of Crelle's journal. Unfortunately, I did not know this article while I was writing my treatise, hence I ommitted to mention it. Recently, you were so kind to send me your article 'Untersuchung der Eigenschaften einer Gattung von unendlichen Reihen' touching this topic. Therefore I would like to thank you and mention a few further results which I obtained in continuation of my work about D.'s series long ago. While being occupied with further questions I was deducted from a rigorous examination of these investigations, and therefore I am not quite sure whether all conditions I suppose on the introduced objects are necessary. The theorems I have found consider the function
$$
(1)\qquad f(s)=\Gamma(s)\left({\Delta\over \pi^2}\right)^{s\over 2}\sum_{n_1,n_2,\ldots n_p} {e^{2\pi i(n_1v_1+n_2v_2+\ldots+n_pv_p)}\over [\varphi(n_1+u_1,n_2+u_2,\ldots n_p+u_p)]^s}
$$
Here $u_1,\ldots u_p,v_1,\ldots v_p$ denote $2p$ real quantities; the summation is with respect to all positive and negative integer values $n_1,n_2,\ldots n_p$, where in case of integers $u_1,u_2,\ldots u_p$ the combination $n_1=-u_1,n_2=-u_2,\ldots n_p=-u_p$ is excluded. Moreover, $\varphi(x_1,x_2,\ldots x_p)$ stands for a positive quadratic form of (without loss of generality) determinant equal to $1$. So far the sum converges for all values of $s$ having real part greater than ${p\over 2}$. The function $f(s)$ can be continued over the whole complex plane and it comes out that this function is univalent and is infinite of the first kind only for $s=0$ and $s=p/2$. If we consider the function
$$
(2)\qquad F(s)=\Gamma(s)\left({\Delta^{p-1}\over \pi^2}\right)^{s\over 2}\sum {e^{-2\pi i(n_1u_1+n_2u_2+\ldots+n_pu_p)}\over [\phi(n_1+v_1,n_2+v_2,\ldots n_p+v_p)]^s}
$$
where $\phi$ is the adjoint form to $\varphi$, then the function $F(s)$ has the same analytic character as $f(s)$ and
\begin{eqnarray*}
f({\textstyle{p\over 2}}-s)&=&e^{-2\pi i(u_1v_1+u_2v_2+\ldots+u_nv_n)}\Delta^{p-1\over 2}F(s)\\ 
f(s)&=&e^{-2\pi i\sum uv}\Delta^{p-1\over 2}F({\textstyle{p\over 2}}-s)\\
F_1(s)&=&\Gamma(s)\left({\Delta^{p(p-2)}\over \pi^2}\right)^{s\over 2}\sum {e^{+2\pi i(n_1v_1+\ldots+n_pv_p)}\over \Delta^{(p-1)s}[\varphi(n_1+u_1,\ldots n_p+u_p)]^s}\\
F_1(s)&=&f(s)\Delta^{(p-2)^2-1}
\end{eqnarray*}
Wishing you pleasant h. I finish with my best reagrds. Your sincerely, A. Hurwitz''\footnote{''Sehr verehrter Herr Geheimrat! Herr Minkowski hat mich schon vor einiger Zeit darauf aufmerksam gemacht, da\ss{} meine Abhandlung \"uber die Dirichlet'schen Reihen in enger Beziehung zu Ihrem im 54 ten Bande des Crelle'schen Journals publizierten Aufsatz steht. Leider war mir dieser Aufsatz zur Zeit als ich jene Abhandlung schrieb nicht bekannt, so da\ss{} ich es unterlie\ss{} auf denselben Bezug zu nehmen. Neuerdings waren Sie so freundlich mir Ihre denselben Gegenstand betreffende Abhandlung 'Untersuchung der Eigenschaften einer Gattung von unendlichen Reihen' zuzusenden. Ich m\"ochte Ihnen hierf\"ur meinen herzlichen Dank aussprechen und mir zugleich erlauben, Ihnen einige weitere Resultate mitzutheilen, welche ich im Anschlu\ss{} an meine Arbeit \"uber die D. Reihen vor l\"angerer Zeit erhalten habe. Durch die Besch\"aftigung mit anderen Fragen bin ich immer wieder von der sorgf\"altigen Durcharbeitung dieser Untersuchungen abgezogen worden, und ich wei\ss{} deshalb nicht, ob die Bedingungen welchen ich die eingef. Gr\"o\ss en unterwerfe s\"amtlich nothwendig sind. Die S\"atze welche ich gefunden habe beziehen sich auf die Function
$$
(1)\qquad f(s)=\Gamma(s)\left({\Delta\over \pi^2}\right)^{s\over 2}\sum_{n_1,n_2,\ldots n_p} {e^{2\pi i(n_1v_1+n_2v_2+\ldots+n_pv_p)}\over [\varphi(n_1+u_1,n_2+u_2,\ldots n_p+u_p)]^s}
$$
Hier bedeuten $u_1,\ldots u_p,v_1,\ldots v_p$ $2p$ reelle Gr\"o\ss en; die Summation bezieht sich auf alle positiven und negativen ganzzahligen Werthe von $n_1,n_2,\ldots n_p$, wobei jedoch falls $u_1,u_2,\ldots u_p$ ganzzahlige Werthe besitzen die Combination $n_1=-u_1,n_2=-u_2,\ldots n_p=-u_p$ auszuschlie\ss en ist. Ferner bedeutet $\varphi(x_1,x_2,\ldots x_p)$ eine positive quadratische Form, deren Determinante (unbeschadet der Allgemeinheit) gleich $1$ angenommen wird. Die Summe convergirt bisher f\"ur alle Werthe von $s$ deren reeller Best. gr\"o\ss er als ${p\over 2}$ ist. Die Function $f(s)$ l\"a\ss t sich \"uber die ganze complexe Zahlenebene fortsetzen und es zeigt sich dabei, da\ss{} diese Function in der ganzen complexen Zahlenebene eindeutig ist und nur f\"ur $s=0$ und $s=p/2$ von der ersten Art unendlich ist. Bildet man nun die Function
$$
(2)\qquad F(s)=\Gamma(s)\left({\Delta^{p-1}\over \pi^2}\right)^{s\over 2}\sum {e^{-2\pi i(n_1u_1+n_2u_2+\ldots+n_pu_p)}\over [\phi(n_1+v_1,n_2+v_2,\ldots n_p+v_p)]^s}
$$
wo $\phi$ die zu $\varphi$ adjungirte Form bezeichnet, so hat die Function $F(s)$ denselben analytischen Charakter wie $f(s)$ und es ist
\begin{eqnarray*}
f({\textstyle{p\over 2}}-s)&=&e^{-2\pi i(u_1v_1+u_2v_2+\ldots+u_nv_n)}\Delta^{p-1\over 2}F(s)\\ 
f(s)&=&e^{-2\pi i\sum uv}\Delta^{p-1\over 2}F({\textstyle{p\over 2}}-s)\\
F_1(s)&=&\Gamma(s)\left({\Delta^{p(p-2)}\over \pi^2}\right)^{s\over 2}\sum {e^{+2\pi i(n_1v_1+\ldots+n_pv_p)}\over \Delta^{(p-1)s}[\varphi(n_1+u_1,\ldots n_p+u_p)]^s}\\
F_1(s)&=&f(s)\Delta^{(p-2)^2-1}
\end{eqnarray*}
Indem ich Ihnen angenehme F. w\"unsche, verbleiben ich mit freundlichen Gr\"u\ss en. Ihr ganz ergebener A. Hurwitz'' Notice that we have not copied scratched formulae or fragments of sentences. The abbreviation 'F.' might stand for 'Ferien' meaning 'holidays', so we have used the letter 'h.' in our translation.}  
\end{quote}
The analytic properties of this zeta-function were proven in Epstein's paper published in 1903 and submitted in January 1912, so 13 years after Hurwitz's letter to Lipschitz. In Epstein's paper \cite{epstein} the notation 
$$
Z\left\vert \begin{array}{ccc} g_1 & \ldots & g_p \\ h_1 & \ldots & h_p \end{array}\right\vert(s)_\varphi:=\sum_{m=-\infty}^{+\infty}{e^{2\pi i\sum_j m_jh_j}\over \varphi((g+m))^{s\over 2}}
$$
is used, where $m=(m_1,\ldots, m_p)$ is a tuple of integers, $g=(g_1,\ldots,g_p)$ and $h=(h_1,\ldots,h_p)$ are tuples of real numbers, and $\varphi((m+g))$ is the quadratic form defined by 
$$
\varphi((x))=\sum_{k,\ell=1}^pc_{k,\ell}x_kx_\ell
$$
with $c_{k,\ell}=c_{\ell,k}$ and $x=(x_1,\ldots,x_p)$. Following Riemann's contour integration method, Epstein derived the following functional equation
\begin{eqnarray*}
\lefteqn{
e^{2\pi i\sum_k g_kh_k}\pi^{-{s\over 2}}\Gamma\left({s\over 2}\right)Z\left\vert \begin{array}{ccc} g_1 & \ldots & g_p \\ h_1 & \ldots & h_p \end{array}\right\vert(s)_\varphi}\\
&=&\Delta^{-{1\over 2}}\pi^{s-p\over 2}\Gamma\left({p-s\over 2}\right)Z\left\vert \begin{array}{ccc} h_1 & \ldots & h_p \\ -g_1 & \ldots & -g_p \end{array}\right\vert(p-s)_\phi. 
\end{eqnarray*}
In the easiest case, when all $h_k$ and $g_k$ are zero and takes the representation $\varphi((m))=m^{\sf t}Qm$ with some symmetric $p\times p$-matrix $Q$ into account, one deduces via substituing $s\mapsto 2s$ the functional equation in the form
$$
\pi^{-s}\Gamma(s)\zeta(s;Q)=(\det Q)^{-1}\pi^{s-{p\over 2}}\Gamma\left({p\over 2}-s\right)\zeta\left({p\over 2}-s,Q^{-1}\right)
$$
for the Epstein zeta-function 
$$
\zeta(s,Q):=\sum_{0\neq m\in\sZ^p} (m^{\sf t}Qm)^{-s}
$$
as it usually appears in the recent literature. We notice that this is essentially in coincidence with Hurwitz's remark. Epstein mentioned Hurwitz's paper \cite{hurwzeta} on the Hurwitz zeta-function as well as Lipschitz \cite{lip89} and Lerch \cite{lerch87} (but, of course, not Hurwitz's diary entry nor his lost letter to Lispchitz). Furthermore, Epstein wrote
\begin{quote}
''Consequently, for the so far most general zeta-function, investigated by Mr Lipschitz, there exists an integral representation in terms of an infinite theta series of {\it arbitrary} characteristic. Following this thought we define a {\it general zeta-function of $p$th order} by a $p$-folded infinite series which generalizes Riemann's zeta-function in just the same way as the general theta series of $p$th order generalizes the elliptic theta series. This extension of the zeta-function is natural since the above mentioned theorem of Riemann about the function $\zeta(s)$ can be transfered to the general zeta-function in a surprisingly easy way.''\footnote{''Dem entsprechend existirt f\"ur die von Herrn Lipschitz untersuchte Function --- die allgemeinste bisher betrachtete 'Zeta'-function --- eine Integraldarstellung mit H\"ulfe einer unendlichen Thetareihe mit {\it beliebiger} Charakteristik. Diesem Gedankengange folgend definiren wir als {\it allgemeine Zetafunction $p^{ter}$ Ordnung} eine $p$-fach unendliche Reihe, welche in gleicher Weise eine Verallgemeinerung der Riemann'schen $\zeta$-function darstellt, wie die allgemeine Thetareihe $p^{ter}$ Ordnung gegen\"uber der elliptischen Thetareihe. Dass diese Erweiterung der Zetafunctionen naturgem\"ass ist, zeigt sich vor allem darin, dass sich der oben angef\"uhrte Riemann'sche Satz \"uber die Function $\zeta(s)$ in \"uberraschend einfacher Weise auf die allgemeinste Zetafunction \"ubertragen l\"asst.'' \cite{epstein}, p. 616} 
\end{quote}
We do not want to go into the details of the appearing notions; and we also do not explain here the different theta-functions but remark that Epstein's approach is obviously via Riemann's second proof of the functional equation (and Vall\'ee Poussin's generalization to Dirichlet $L$-functions).\footnote{In a subsequent paper the authors intend to investigate the Hurwitz estate further in order to find details about his approach to zeta-functions associated with quadratic forms.}
\par

A last comment about Hurwitz's letter to Lipschitz. In 1887, after studying with Hurwitz in K\"onigsberg, his former pupil and friend Hermann Minkowski became professor at the University of Bonn where Lipschitz was professor (from 1864 until the end of his life in 1903). We may speculate that this contact had been a reason for the correspondence between Hurwitz and Lipschitz.\footnote{In 1897, Minkowski married Auguste Adler from Strasbourg, the place where Epstein obtained his docotrate in 1895 (advised by Elwin Christoffel) and remained as docent until 1919 when he received the call of the recently founded university in his native town of Frankfurt/Main. As a matter of fact, Minkowski was frequently visiting his brother and famous physician Oskar at Strasbourg starting from 1889 as we know from Hermann Minkowski's letters to Hilbert \cite{zass}; during some of these visits Minkowski got in touch with the mathematicians at Strasbourg University and Christoffel in particular (see \cite{zass}, p. 34, 36). A word about Frankfurt University, too. This institution was founded in 1914, mainly by donations of the Jewish community of Frankfurt, but the Jew Epstein was driven to commit suicide under the Nazi regime. Actually, all members of the Frankfurt Mathematical Seminar except Siegel lost their positions. For details about Epstein's life we refer to Siegel \cite{siegelh} and Burde et al. \cite{burdeh}.} Bianchi's hint from his 1881 letter is not mentioned.  
\smallskip

We conclude with a discussion about the name attributed to Hurwitz's generalization. Reading Landau's historical survey \cite{landau} from 1906, we find the Hurwitz zeta-function with notation $\zeta(s,w)$; a quotation of Hurwitz is missing, however, Mellin \cite{mellin1} is mentioned. In 1929, Walter Schnee \cite{schnee} gave a new proof of the functional equation for $\zeta(s,\alpha)$ and Dirichlet series with periodic coefficients, however, Hurwitz is not mentioned. The same in Helmut Hasse's paper \cite{hasse}, where $\zeta(s,\alpha)$ is called 'allgemeine $\zeta$-Reihe'. In 1950/51 this was changed by two articles in the same volume of the Proceedings of the American Mathematical Society, namely Nathan J. Fine \cite{fine} and Tom M. Apostol \cite{apostol}. The first paper, the one by Fine, gives a proof of the functional equation for the Hurwitz zeta-function by means of Riemann's second approach using Poisson summation and properties of Jacobi's theta-function. The second paper, the one by Apostol, gives another proof relying on Lerch's transformation formula. Both authors explicitly call $\zeta(s,\alpha)$ the Hurwitz zeta-function. We do not consider the many further research papers on the Hurwitz zeta-function here. 
The name Epstein zeta-function for the generalization associated with quadratic forms, however, has been used since Edward Charles Titchmarsh's article \cite{titch} from 1934.
\par

The Hurwitz zeta-function $\zeta(s,\alpha)$ with an irrational parameter $\alpha$ is not treated in Hurwitz's original paper. For example, in the monographies of Lang \cite{lang0,lang} and Garunk\v stis \& Laurin\v cikas \cite{garunklauri} the reader may find (essentially the same) proof for the Hurwitz zeta-function with a not necessarily rational parameter and in the latter one the proof of the functional equation for the more general Lerch zeta-function.
\smallskip

Obviously, some of the mathematicians mentioned above were not very well informed about what others did. The standards of citation were rather different in the 19th century, and so the distribution of research results and of journals as well. Another aspect is that it was not quite clear in the process of generalization from the Riemann zeta-function and Dirichlet $L$-series which one of all the various possible generalized zeta-functions would be the best (if one really would like to coin the notion 'best' in this context). The struggle for the {\it best} generalization of the Riemann zeta-function reminds the authors on the numerous approaches to multiple (multivariate) zeta-functions and zeta-values in the literature (e.g., the Euler -Zagier multiple zeta-functions) investigated intensively for about 25 years. 

\section{K\"onigsberg \& Zurich --- Entire Functions}

In 1884 Adolf Hurwitz became extraordinary professor at the University of K\"onigsberg (now Kaliningrad in Russia) on invitation of Ferdinand Lindemann. Although friends and relatives as well as most of his former mathematical friends and colleagues were far distant Hurwitz spent a very good time in K\"onigsberg. He produced many beautiful theorems and he was rather successful in teaching: among the few K\"onigsberg students Hurwitz had with David Hilbert and Hermann Minkowski two brilliant pupils. During frequent walks Hurwitz introduced the two students to various mathematical disciplines and guided their first steps in research. The three of them became lifelong friends. Moreover, Adolf Hurwitz got to know his later wife, Ida Samuel, the daughter of the professor for Pathology at K\"onigsberg University; they married in summer 1892. In the same year Frobenius left Zurich for Berlin and Hurwitz was offered the vacant ordinary professorship. The young family moved to Switzerland.

\subsection{Prehistory: The Proof of the Prime Number Theorem}

Riemann's celebrated paper \cite{rie} is at least nowadays considered most important with respect to its not rigorously proved statements and conjectures. This memoir has indeed been considered as a signpost for developing complex analytic methods in order to reveal distribution properties of the prime numbers. Taking Riemann's personality into account one may understand his contribution also as an outline of a far reaching analytic method in number theory, leaving the details to the posterity. The announced product representation was obtained by Jacques Hadamard (1865-1963) \cite{hada93} in 1893, the explicit formula linking the zeros of $\zeta(s)$ with the prime numbers was proved by von Mangoldt \cite{mango95} in 1895, the so-called Riemann-von Mangoldt formula for the number of zeros inside certain rectangles was as well obtained by von Mangoldt in the same paper; only the famous Riemann hypothesis claiming that all nontrivial zeros lie on the critical line has remained unsolved until today.\footnote{Actually, Riemann wrote that it is 'very likely' that all these zeros lie on a line; he did neither announce this comparable to his other claims nor did he formulate it as a conjecture. Nowadays we can only speculate whether Riemann knew about the arithmetical consequences of such a zero distribution.} Building on Riemann's ideas, Hadamard \cite{hada96} and Charles de la Vall\'ee Poussin (1866-1962) \cite{poussin96} solved in 1896 independently a conjecture of Gauss on the number of primes below a given magnitude; their reasoning culminated in the celebrated prime number theorem giving the asymptotics for the counting function $\pi(x)$ of prime numbers $p\leq x$, namely 
$$
\lim_{x\to\infty}\pi(x){\log x\over x}=1\, ,\qquad\mbox{resp.}\qquad \pi(x)\sim {x\over \log x}
$$
as $x\to \infty$. More precise formulae provide explicit error terms depending on the zero-free region for the zeta-function. For details on the history of the development of the prime number theory we refer to Narkiewicz's monumental \cite{narkie}, Laugwitz's Riemann biography \cite{laugw}, and Schwarz's survey \cite{wschwarz}.


\subsection{Hadamard's Theory of Entire Functions of Finite Order} A function $f(s)$ is called entire if it is analytic (holomorphic) throughout the whole complex plane; if its absolute value is not growing too quickly as $\vert s\vert$ tends to infinity, it is said to be of finite order (and 'of finite genre' in older literature). Since the Riemann zeta-function $\zeta(s)$ is analytic except for a simple pole at $s=1$, it is natural to consider $(s-1)\zeta(s)$ or similar functions with respect to the general theory of entire functions. This line of inquiry had been established among others by Jacques Hadamard in his dissertation \cite{hada92} from 1892 advised by \'Emile Picard. Hadamard's motivation, however, was not the theory itself but its application to the unsolved problems around the Riemann zeta-function. In fact, his theory of entire functions of finite order constitutes a major ingredient in his and de la Vall\'ee Poussin's celebrated proof of the prime number theorem and later refinements. 
\par

In the {\it Jahrbuch der Mathematik} we can read in a review by Adolf Hurwitz concerning Hadamard's article \cite{hada93}:
\begin{quote}
''The article under review is the accomplishment of a work which had been awarded the great mathematical prize of the Parisian Academy of Sciences. Taking into account that the author had had only a limited amount of time, it is attributed that the formulation of the article causes trouble to its comprehension. For example, the proof of the main theorem from the first part should have been arranged differently in order to be rigoros. This is just about the outward appearance of the work; as regards content it may be considered as one of the most important function-theoretical pieces of the last years.''\footnote{''Die vorliegende Abhandlung ist die weitere Ausf\"uhrung einer Arbeit, die von der Akademie der Wissenschaften zu Paris mit dem grossen mathematischen Preise gekr\"ont worden ist. Dem Umstande, dass der Verfasser nur eine beschr\"ankte Zeit zur Verf\"ugung hatte, ist es wohl zuzuschreiben, dass die Abfassung der Abhandlung dem Verst\"andnisse manche Schwierigkeiten bereitet. So muss beispielsweise der Beweis des Hauptsatzes im ersten Teile der Arbeit anders angeordnet werden, um bindende Kraft zu erhalten. Diese Ausstellung betrifft indessen nur die \"aussere Form der Arbeit; inhaltlich darf dieselbe wohl als eine der bedeutendsten functionentheoretischen Arbeiten der letzten Jahre bezeichnet werden.''} 
\end{quote}
Indeed, Hadamard's article under consideration had been awarded the {\it Grand Prix des Sciences Math\'ematique} of the French Academy of Sciences in 1892 for a treatise on the Riemann zeta-function (reporting about the results from his doctorate). Hurwitz wrote several hundreds of reviews of articles in German, French, Italian, and English for the {\it Jahrbuch \"uber die Fortschritte der Mathematik}\footnote{which was founded as early as 1868 by Carl Orthmann and Felix M\"uller and served as pre-runner of Mathematical Reviews and Zentralblatt until 1942; now its may reviews are included into the Zentralblatt data bank.} and these reviews provide a mathematical fingerprint of his research activities and interests. Hurwitz followed the many works dealing with entire functions in particular and his diaries contain numerous entries about this topic. In his Mathematical Diary No. 15, we find for 22 February 1897 an entry analysing the Hadamard's papers \cite{hada962,hada96} on the distribution of zeros of $\zeta(s)$ and its arithmetical consequences. It is interesting to notice that Hadamard cited Hurwitz's paper \cite{hurwzeta} (besides Cahen \cite{cahen2}) and made use of his identity. 

\subsection{P\'olya and Hurwitz's Estate}

George (Gy\"orgy) P\'olya was born in Budapest 1887; he finished his studies 1912 with a doctorate supervised by Lip\'ot Fej\'er in Budapest, investigating geometric probability. After a postdoc in G\"ottingen he started at the Eidgen\"ossische Hochschule Zurich in 1914 on promotion by Hurwitz. He became extraordinary professor in 1920 and ordinary professor only in 1928 (despite his numerous important results).\footnote{The Lerch pupil Michel Plancherel was chosen to inherit Hurwitz's chair, not P\'olya although he was considered as outstanding mathematician by none less than Hilbert and Hadamard; see Frei \& Stammbach \cite{stamm}, p. 51, for details.} In 1940 P\'olya moved to the U.S.A. where he had professorships at several universities, since 1942 in Stanford. P\'olya died in Palo Alto in 1985. He is considered as one of the outstanding mathematicians of the twenteeth century. As his mentor and colleague Hurwitz P\'olya had a broad mathematical knowledge and interest.
\begin{quote}
''He also credited Hurwitz's extensive mathematical diaries for a number of problems that later appeared in the {\it Aufgaben und Lehrs\"atze}. These diaries may have influenced P\'olya in other ways: he kept a mathematical log himself throughout much of his career, recording mathematical conversations and keeping track of his own mathematical ideas.'' 
\end{quote}
wrote Alexanderson \& Lange \cite{alex}. The mentioned book {\it Aufgaben und Lehrs\"atze aus der Analysis} was joint work with Gabor Szeg\H{o} and the innovative style of this collection of exercises (ordered with respect to methods to prove these results rather than the results itself) is still a pearl in the mathematical literature. Besides many interesting tasks it includes, for example, the proof of the existence of infinitely many prime numbers by use of the coprimality of the Fermat numbers $2^{2^n}+1$ which is probably based on lectures of Hurwitz' which were later published as \cite{hurwitzbuch}\footnote{and may have been also an inspiration for the book \cite{posz} of P\'olya \& Szeg\"o} and the letter exchange between Goldbach and Euler; for details we refer to Haas \cite{haas}. 
\par 

Starting in 1921, two years after Adolf Hurwitz's death in Zurich 18 November 1919, P\'olya became literary executor of Hurwitz's estate. In 1933 Hurwitz's collected papers were finally published, edited by P\'olya. The mathematical diaries of Adolf Hurwitz definitely had been another source of inspiration for P\'olya.\footnote{Since in Hurwitz's last diary one can find some notes from 1918 about Arthur Cayley's studies of counting trees, namely alkane ${\sc C_nH_{2n+2}}$ with certain restrictions, one could speculate whether this might have been an inspiration for P\'olyas celebrated enumeration theory.} Here we shall focus on zeta-function theory.
\begin{quote}
''Of those topics which are treated in the diaries for several years again and again we shall mention here Riemann's zeta-function as an example; in the published works of Hurwitz this object is not mentioned at all. Very early Hurwitz had tried a path towards the Riemann hypothesis in vain which some years later also Jensen found: (...) The aim was unreachable.''\footnote{''Von denjenigen Gegenst\"anden, auf welche die Tageb\"ucher lange Jahre hindurch immer wieder und wieder zur\"uckkommen, sei hier als Beispiel nur die Riemann'sche $\zeta$-Funktion hervorgehoben; in den ver\"offentlichten Arbeiten von Hurwitz wird dieser Gegenstand wohl nicht einmal erw\"ahnt. Schon fr\"uh hat Hurwitz den Weg zur Riemann'schen Vermutung versucht, den einige Jahre sp\"ater auch Jensen, ebenfalls ergebnislos versucht hat: (...) Das Ziel war unerreichbar.''}
\end{quote}
These words are due to George P\'olya introducing and commenting the unpublished notes of Adolf Hurwitz \cite{diaries}, p.6. The first mathematical diary starts 25 April 1882, too late for Hurwitz's work on the Hurwitz zeta-function \cite{hurwzeta}. 


\begin{quote}
''Hurwitz was the type of mathematician who always achieved something when seriously using his power, if not the envisaged objective, then at least something interesting besides. Two side issues of his efforts about the $\zeta$-function reached the mathematical community by oral communication:...''\footnote{''Nun war Hurwitz von dem Mathematikertypus, der bei ernstlicher Einsetzung seiner Kr\"afte immer etwas erreicht, wenn nicht das urspr\"unglich ins Auge gefasste Ziel, so doch etwas Interessantes auf den Seitenwegen. Zwei Nebenresultate seiner Bem\"uhungen um die $\zeta$-Funktion sind durch m\"undliche Mitteilungen in die mathematische Oeffentlichkeit gelangt: ...''}
\end{quote}
This quotation is again from P\'olya's introductary words to Hurwitz's estate \cite{diaries}, p.7. 
Indeed P\'olya himself had continued some of Hurwitz's ideas. 
\par

In 1914, Godfrey Harold Hardy \cite{hardy} proved that infinitely many zeros of $\zeta(s)$ lie on the critical line. His reasoning was based on the real-valued function defined essentially by the left hand-side of the functional equation, i.e.,
$$
\Xi(t):={\textstyle{1\over 2}}s(s-1)\pi^{-{s\over 2}}\Gamma({\textstyle{s\over 2}})\zeta(s)\LARGE\vert_{s={1\over 2}+it}
$$
and its integral representation
$$
\Xi(t)={\textstyle{1\over 2}}-(t^2+{\textstyle{1\over 4}})\int_1^\infty x^{-{3\over 4}}\cos({\textstyle{1\over 2}}t\log x)\sum_{n\geq 1}\exp(-\pi n^2 x)\d x.
$$
Hardy's reasoning relies on the integral
$$
\int_0^\infty \Xi(t){t^{2n}\over t^2+{1\over 4}}\cosh(\alpha t)\d t 
$$
and taking the limit $\alpha\to {\pi\over 4}-$ under the assumption that $\Xi(t)$ is of constant sign for all sufficiently large $t$. By the upcoming contradictions it follows that $\Xi(t)$ must have infinitely many odd order zeros, and since all factors of $\Xi(t)$ do not vanish except $\zeta({1\over 2}+it)$ possibly, this proves the existence of infinitely many zeta zeros on the critical line. In Hurwitz's Mathematical Diary no. 26 we can find an outline of Hardy's proof.
P\'olya \cite{polya27} found a variant of Hardy's proof by the associated Fourier integral
$$
\Phi(u):={\textstyle{1\over \pi}}\int_0^\infty \Xi(t)\cos(ut)\d t 
$$
and its Taylor expansion. The same line of inquiry has been treated by P\'olya in his papers \cite{polya26a,polya26b,polya27b}. Probably, the easiest proof of Hardy's theorem is due to Landau \cite{landau27}, pp. 78, which relies on the inequality
$$
\int_T^{2T}Z(t)\d t<\int_T^{2T}\vert Z(t)\vert\d t,
$$ 
valid for sufficiently large $T$, where 
$$
Z(t):=\exp(i\theta (t))\zeta({\textstyle{\frac{1}{2}}}+it)
$$
with 
$$
\theta (t):=\pi ^{-{\frac{it}{2}}}{\Gamma ({\textstyle{\frac{1}{4}}}+{\textstyle{\frac{it}{2}}})\over|\Gamma ({\textstyle{\frac{1}{4}}}+{\textstyle{\frac{it}{2}}})|}
$$
is real-valued. Hardy's approach allows quantitative results, however, the so far best quantitative bound for the number of zeros on the critical line was obtained by Levinson's method and its technical refinements.\footnote{see Iwaniec \cite{iwaniec} for a recent presentation}
\par

P\'olya's paper \cite{polya18} from 1918 deals with the zeros of entire functions represented by certain Fourier integrals. It even contains a theorem due to Hurwitz and his proof (see \cite{polya18}, p. 368), namely: {\it given an even function $f(t)$ defined for $t\in(-1,1)$ having alternating Fourier coefficients, i.e.,
$$
f(t)\sim {\textstyle{1\over 2}}a_0+a_1\cos(\pi t)+a_2\cos(2\pi t)+a_3\cos(3\pi t)+\ldots
$$
with $a_0>0, a_1<0, a_2>0, a_3<0,\ldots$, then the zeros of the entire function 
\begin{equation}\label{uz}
U(z)=\int_0^1f(t)\cos(zt)\d t
\end{equation}
are all real and simple, and such that in each of the intervals $\ldots, (-2\pi,-\pi), (-\pi,0),(0,\pi),(\pi,2\pi),\ldots$ contains exactly one zero.} This result and its rather lengthy proof is related to P\'olya's reasoning; moreover, a paper of Jensen is mentioned. And the name Jensen comes up in another context as well. For his article \cite{polya27} P\'olya has studied the mathematical papers left by Johan Ludwig William Valdemar Jensen (1859-1925).\footnote{It is interesting to notice that around the same time Carl Ludwig Siegel was investigating Riemann's estate in G\"ottingen; his study \cite{siegel} has changed the impression many contemporaries had about Riemann's approach to number theory.} The autodidact Jensen never had an academic appointment, however, working for the International Bell Telephone Company at Copenhagen, he found time for severe mathematical research. In zeta-function theory the so-called Jensen formula plays an important role in estimating the number of zeros of an entire function in a given disk.\footnote{According to Narkiewicz \cite{narkie}, p. 258, Jensen's formula was essentially known by Carl Gustav Jacobi already in 1827.} In 1911, during the Second Congress of Scandinavian Mathematicians at Copenhagen Jensen \cite{jv} announced five articles about new algebraic and function-theoretical methods with applications to the Riemann zeta-function (JV), however, none of these announcements was published until Jensen's death in 1925. Probably inspired by his studies of Hurwitz's estate, P\'olya investigated the Jensen estate (JN) but could not find what he was aiming at.
\begin{quote}
''This all is to find in several small notebooks and a great number of loose pages such that in most cases it is impossible to fix in which sequence and which time the notes have been written down. As follows from this description the state of the estate does not allow to set its content with certainty.''\footnote{''Dies alles in einigen kleinen Heften und an einer grossen Anzahl von losen Bl\"attern, so dass in den allermeisten F\"allen unm\"oglich festzustellen ist, in welcher Reihenfolge und zu welchem Datum die Aufzeichnungen entstanden sind. Wie aus dieser Schilderung hervorgeht, erlaubt der Zustand des Nachlasses nicht, dessen Inhalt mit Sicherheit festzustellen.''} 
\end{quote}
It seems that Jensen was in some sense a messy mathematician as Hurwitz was a tidy character. This chaos did not stop P\'olya's investigations but led him to write up an at that time complete survey what was known about integral representations of the Riemann zeta-function. In this context P\'olya gave credit to Hurwitz:
\begin{quote}
''I is contained in JV (p. 188), II and III are often quoted in JN. For II and IV see $P_3$. --- II and III were (as well as the consequences mentioned in Nr. 6 following Theorem III) known by several friendly mathematicians; both appear in the mathematical diaries of A. Hurwitz (preserved by the library of the Eidg. Technischen Hochschule in Zurich) with date 1899.''
\footnote{''I ist in JV (S. 188), II und III sind in JN \"ofters hervorgehoben. F\"ur II und IV vgl. $P_3$. --- II und III waren (nebst der in Nr. 6 im Anschluss an Satz III hervorgehobenen Konsequenzen) mehreren mit mir befreundeten Mathematikern bekannt; beide kommen in den Tageb\"uchern von A. Hurwitz (aufbewahrt in der Bibliothek der Eidg. Technischen Hochschule in Z\"urich) unter dem Datum 1899 vor.'' \cite{polya27}, p. 11. Here $P_3$ is an abbreviation for \cite{polya26a}.}
\end{quote}
Indeed, in the Mathematical Diary no. 17, p. 112-115, we can find the corresponding formula (see Figure 2).
\begin{figure}[h]
\includegraphics[height=11.2cm]{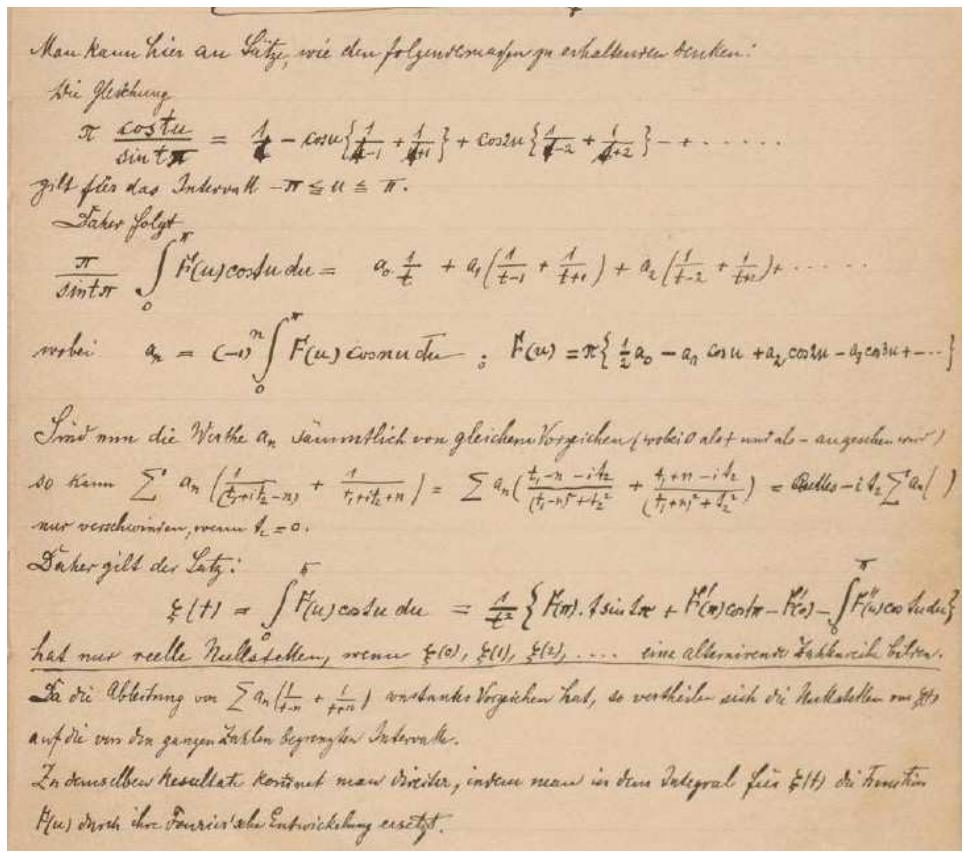}
\caption{Hurwitz's Diary no. 17, p. 115. Formula (\ref{uz}) appears in disguised form in the lower third. Another entry of the diary, a few pages later, bears the date 3 December 1899.}
\end{figure}

Now we are in the position to explain what P\'olya meant in writing ''Very early Hurwitz had tried a path towards the Riemann hypothesis in vain which some years later also Jensen found'' (as mentioned above). In 1889, Hurwitz \cite{hurwitz90} extended a classical result due to Sim\'eon Poisson \cite{poisson}, p. 178, on the zeros of Bessel functions $J_n$ (whose definition we omit here), namely that all the zeros are real when $n\geq 1$. He extended and quantified Poisson's result about the zero-distribution using Rouch\'e's theorem\footnote{which probably was unknown to him and better should be named after Augustin Louis Cauchy as already Eug\`ene Rouch\'e had remarked; see Bottazini \& Gray \cite{bg} for details.} and more advanced techniques; in particular he obtained his famous theorem on {\it uniformly convergent sequences of meromorphic functions having limit either identically vanishing or}, more interestingly, {\it every zero of the limit function is limit of the zeros of the functions in the sequence}. This result is now a standard tool in the complex analytic tool box. Using this observation, Hurwitz was led to a separation of the real zeros as in the case of (\ref{uz}) above. And here is the link to his approach to the Riemann hypothesis. Hurwitz wrote:
\begin{quote}
''The proof of this theorem is, following Poisson's process, usually based on an integral representation...''\footnote{''Der Beweis dieses Satzes pflegt man nach Poisson's Vorgang auf eine Integralformel zu gr\"unden...''; \cite{hurwitz90}, p. 246.}
\end{quote}
We observe the similarity with the situation considered by Hardy and the entry in Diary no. 17 (Figure 2) from above. In order to prove that all nontrivial zeros of the zeta-function are on the critical line one could aim at finding a suitable integral representation; this could be realized by giving criteria that an entire function (e.g., a Fourier transform of a kernel) has all its zeros on the real axis. This idea appears in the literature several times; however, it seems that Hurwitz was among the first to have considered this approach ''in vain (...) the aim was unreachable'' as P\'olya remarked.
\par

In a subsequent paper \cite{polya27b}, following an idea of Landau, P\'olya obtained the related remarkable result that all zeros of the function 
$$
\Xi^*(t):=4\pi^2\int_0^\infty \cosh({\textstyle{9\over 2}}u)\exp(-2\pi\cosh (2u))\cos(ut)\d t
$$
are real; although this function is an approximation of $\Xi(t)$, unfortunately this has no impact on the location of the zeŕos of the Riemann zeta-function. Nevertheless it reminds us of another related approach to the Riemann hypothesis. Based on a theorem due to Edmond Nicolas Laguerre on entire functions with real zeros, P\'olya \cite{polya13} studied sequences of polynomials with positive real roots and proved that a uniform limit is either vanishing identically or an entire function of order zero or one (see (\ref{genus}) below). Interestingly, P\'olya's paper includes a certain criticism of Laguerre's proof:
\begin{quote}
''The correctness of this proof can be doubted; I am of the opinion that the proof can be realized rigorously without any substantial changes. --- The here given proof is completely independent of the correctness of Laguerre's reasoning.''\footnote{''Die Richtigkeit dieses Beweises kann angezweifelt werden; ich bin der Ansicht, dass der Beweis ohne wesentliche Aenderungen des Gedankenganges l\"uckenlos zu f\"uhren ist. --- Der hier gegebene Beweis ist von der Richtigkeit des Laguerreschen v\"ollig unabh\"angig, und beruht auf einer ganz anderen Grundlage.'', \cite{polya13}, p. 280} 
\end{quote}
It should be mentioned that \cite{polya13} had been written in 1913, P\'olya's last year in G\"ottingen before he moved to Zurich. Around the same time Jakob Grommer wrote his doctorate \cite{grommer} supervized by Hilbert on this topic. He found an improvement of one of Hurwitz's unpublished theorems on entire functions with real zeros only; the original result and questions about its possible generalization were communicated by Otto Toeplitz and it was Hilbert who supported Grommer ''who apart from his physical ailment faced triple jeopardy as a foreigner, orthodox Jew, and as someone without the equivalent of the {\it Abitur}, since he graduaded from a Talmudic school rather than a traditional Gymnasium'' wrote Rowe \& Schulmann \cite{rowe}. Grommer's thesis forms the second issue P\'olya had mentioned in his comment about the Hurwitz estate.  
\par

It seems that these kind of questions were among the first that Hurwitz and P\'olya discussed after the latter's arrival at Zurich in 1913. In Diary no. 26 one can find an entry on 'P\'olya's theorem' under date 26 February 1914 which is related to his publication \cite{polya13a}; actually, the same diary contains a note with respect to Laguerre's theorem as well as an application of Charles Biehler's theorem \cite{biehler} on real zeros of entire functions with comment 
\begin{quote}
''A pretty application of Biehler's theorem came to me last night.''\footnote{''Eine h\"ubsche Anwendung des Biehler'schen Satzes fiel mir vergangene Nacht ein.''}
\end{quote}
\par

In the recent literature the so-called Laguerre--P\'olya class is defined as the set of entire functions $E(s)$ being locally limits of a series of polynomials having only real roots. Such functions obey (like polynomials) the Hadamard product representation
\begin{equation}\label{genus}
E(s)=s^d\exp(a+bs+cs^2)\prod_n\left(1-{s\over s_n}\right)\exp\left({s\over s_n}\right),
\end{equation}
where $a,b,c$ are constant, $b$ and $c$ real, the latter one zero or negative, $d$ is the multiplicity of a (hypothetical) zero at $s=0$, and the $s_n$ denote the further zeros of $E$. There is a lot of literature on this class of entire functions and their zero-distribution; the impact on the Riemann zeta-function and allied functions has been investigated for decades, e.g.,  
Louis de Branges \cite{branges} should be mentioned here. However, since this is not really related to Adolf Hurwitz and his mathematics, we refer here only to the survey article \cite{csordas} by Csordas, Norfolk \& Varga. In particular P\'olya's paper \cite{polya13} has been extended and precised by Dennis Hejhal \cite{hejhal}. 
\par

Another idea in the context of the Riemann hypothesis related to P\'olya is the so-called Hilbert-P\'olya conjecture that the ordinates of the nontrivial zeros of the Riemann zeta-function correspond to eigenvalues of an unbounded self-adjoint operator. This idea as well does not go back to conversations with Hurwitz but arose from P\'olya's time in G\"ottingen and a question by Landau.\footnote{See the correspondence \cite{odlyzkopolya} between Andrew Odlyzko and P\'olya on this topic.} 

\section{Concluding Remarks}

Introducing analytic methods to number theory was a change of paradigm. The ideas of Euler, Dirichlet, and Riemann, however, were ahead the methods at their time. Open question about prime number distribution became a driving force for the development of certain aspects of complex analysis. The distribution of zeros of zeta- and $L$-functions is still not understood very well. This might have two reasons: it is neither known what impact a Riemann type-functional equation on the zeros distribution actually has nor is it clear why arithmetically relevant Dirichlet series with an Euler product seem to have all its complex zeros on one single line.\footnote{Actually, the notion of an $L$-function is usually used for zeta-functions having an Euler product. An axiomatic setting for such arithmetically relevant Dirichlet series was provided by Atle Selberg \cite{selberg}. It is legend that the letter '$L$' is chosen with respect to Dirichlet's given name '$L$eujene'.} The functional equation alone cannot be the driving force for the complex zeros to lie on the critical line. Harold Davenport \& Hans Heilbronn \cite{davenport1} showed that the Hurwitz zeta-function $\zeta(s,\alpha)$ has infinitely many zeros in the half-plane $\Re s>1$ provided $\alpha\neq {1\over 2},1$ is rational or transcendental; moreover, they showed the explicit existence of Dirichlet series satisfying a Riemann-type functional equation having infinitely many zeros off the critical line. It is expected that the Euler product should have a severe impact on the zero-distribution (although there is no absolute convergence of the Euler product inside the critical strip).  
\par 
 
The second proof Riemann gave for the functional equation of the zeta-function may be considered more important by its relation to modular forms. This line of investigation has been studied by Erich Hecke starting with his proof of the functional equation for Dedekind zeta-functions \cite{hecke0} in 1916/17 and the relation between Dirichlet series having a functional equation of the Riemann type and modular forms \cite{hecke} in 1936. Nowadays the role of Dirichlet series for arithmetical questions is quite fairly understood. We quote 
\begin{quote}
''In recent years $L$-functions and their analytic properties have assumed a central role in number theory and automorphic forms.''
\end{quote}
from the abstract of Gelbart \& Miller \cite{gelb} which gives an excellent reading for the modern zeta-function theory. We shall briefly mention here John Tate's thesis \cite{tate} in which a proof of the functional equation is given on the basis of local methods. This method can be interpreted as an analytic variation of the local global principle discovered by Helmut Hasse for quadratic forms. Another line of inquiry has been suggested by Robert P. Langlands in the 1970s with his 'Langlands Program' claiming that all arithmetically relevant $L$-functions (automorphic or motivic) arise as products of $L$-functions associated with automorphic representations; the celebrated proof of the Shimura \& Taniyama conjecture due to Andrew Wiles et al. is a prominent example. For details we refer once again to \cite{gelb}. On the contrary, the first proof Riemann gave in his path-breaking paper and Hurwitz's generalization still gives (at least in the authors' opinion) the most simple approach to Dirichlet $L$-functions (without using Poisson's summation formula).
\par

We summarize. Shortly after receiving his doctorate in 1881 Adolf Hurwitz made new contacts and found research problems in new mathematical topics. One of them, zeta-function theory, should be a field of interest for the rest of his life although this is apart from his paper on the Hurwitz zeta-function not visible in his list of publications. His mathematical diaries and a few reviews as well as the remarks in P\'olya's papers, however, give proof for this liking. Hurwitz's ideas concerning the Riemann zeta-function and its relatives have inspired many contemporaries and later mathematicians and traces can be found until today. Last not least, Hurwitz was in possession of fundamental insights about an important generalization of the variety of zeta-functions under investigation at the turn of the century. Unfortunately, he (probably) only communicated his results about these zeta-function to quadratic forms in a letter to Lipschitz but did not find time to polish and publish his notes; nowadays all the credit is given to Epstein. 
\par

We conclude with another quotation from P\'olya's text introducing the mathematical diaries of Adolf Hurwitz:
\begin{quote}
''Perhaps one will find in the diaries some more isolated pretty results and promising germ...''\footnote{''Wohl wird noch in den Tagb\"uchern manches h\"ubsche Einzelresultat, mancher entwicklungsf\"ahige Anregungskeim zu finden sein...''; \cite{diaries}, p. 7}
\end{quote}
\medskip

\noindent {\bf Acknowledgements.} The photographs were taken from internet sources, and some of them are courtesy of the ETH library; most of the pictures stem from the webpages of the MacTutor History of Mathematics Archive at St Andrews University, Scotland. We are grateful to these institutions for their support. Moreover, we would like to thank Klaus Volkert and J\"urgen Wolfart for reading a former version of this article and for giving valuable hints. Last not least, we would like to express our gratitude to the organizers J\"urgen Sander, Martin Kreh, and Jan-Hendrik de Wiljes of the conference {\it Elementare und analytische Zahlentheorie} at Hildesheim University in 2014 for the kind hospitality and their support to publish this volume in memory of Wolfgang Schwarz.

\bigskip

\noindent {\footnotesize Nicola Oswald, Department of Mathematics and Informatics, University of Wuppertal, Gau\ss str. 20, 42\,119 Wuppertal, Germany, oswald@uni-wuppertal.de;\\
\quad and\\
Department of Mathematics, W\"urzburg University, Emil-Fischer-Str. 40, 97\,074 W\"urzburg, Germany, nicola.oswald@mathematik.uni-wuerzburg.de}

\noindent {\footnotesize J\"orn Steuding, Department of Mathematics, W\"urzburg University, Emil-Fischer-Str. 40, 97\,074 W\"urzburg, Germany, steuding@mathematik.uni-wuerzburg.de}

%
%
%


\small

\begin{thebibliography}{9}

\bibitem{alex}{\sc G.L. Alexanderson, L.H. Lange}, Obituary George P\'olya, {\it Bull. London Math. Soc.} {\bf 19} (1987), 559-608

\bibitem{apostol}{\sc T.M. Apostol}, Remark on the Hurwitz Zeta Function, {\it Proc. A.M.S.} {\bf 2} (1951), 690-693

\bibitem{bianchibrief}{\sc L. Bianchi}, Letters of Luigi Bianchi to Adolf Hurwitz, Korrespondenz von Adolf Hurwitz, Cod Ms Math Arch 76, Nieders\"achsische Staats- und Universit\"atsbibliothek G\"ottingen 

\bibitem{bianchi}{\sc L. Bianchi}, {\it Opere, Vol. XI: Corrispondenza}, Edizioni Cremonese, Roma 1959 

\bibitem{biehler}{\sc Ch. Biehler}, Sur une classe d'equation alg\'ebraiques dont toutes les racines sont r\'eelles, {\it Nouvelles annales math\'ematique} {\bf 19} (1880), 149-152

\bibitem{humbo}{\sc K.-R. Biermann} (ed.), {\it Briefwechsel zwischen Alexander von Humboldt und Peter Gustav Lejeune Dirichlet}, Akademie der Wissenschaften der DDR, Berlin 1982

\bibitem{bg}{\sc U. Bottazini, J. Gray}, {\it Hidden Harmony -- Geometric Dantasies}, Springer 2013

\bibitem{branges}{\sc L. de Branges}, {\it Hilbert spaces of entire functions}, Prentice-Hall, London 1968

\bibitem{burdeh}{\sc G. Burde, W. Schwarz, J. Wolfart}, Max Dehn und das Mathematische Seminar der Universit\"at Frankfurt im Fluss der Zeit, {\it Algorismus} {\bf 44} (2004), 462-483

\bibitem{cahen2}{\sc E. Cahen}, Sur la fonction $\zeta(s)$ de Riemann et sur des fonctions analogues, {\it Ann. de l'\'Ec. Norm.} {\bf XI} (1894), 75-164

\bibitem{cantor}{\sc M. Cantor}, Nachruf an Oskar Schl\"omilch, {\it Bibl. math.} {\bf 3} (1901), 260-281 

\bibitem{clausen}{\sc T. Clausen}, Beweis des Schl\"omilch Lehrsatzes, {\it Arch. Math. Phys.} {\bf 30} (1858), 166-170

\bibitem{csordas}{\sc G. Csordas, T.S. Norfolk, R.S. Varga}, The Riemann hypothesis and the Tur\'an inequalities, {\it Trans. Amer. Math. Soc.} {\bf 296} (1986), 521-541

\bibitem{daven}{\sc H. Davenport}, {\it Multiplicative Number Theory}, Springer 1967

\bibitem{davenport1}{\sc H. Davenport, H. Heilbronn}, On the zeros of certain Dirichlet series I,II, {\it J. London Math. Soc.} {\bf 11} (1936), 181-185; 307-312

\bibitem{dirichlet37}{\sc P.G.L. Dirichlet}, Beweis des Satzes, dass jede unbegrenzte  arithmetische Progression, deren erstes Glied und Differenz ganze Zahlen ohne gemeinschaftlichen Factor sind, unendlich viele Primzahlen enth\"alt, {\it Abhandlungen Kgl. Preu\ss . Akad. Wiss.} (1837), 45-81; Werke I, Reimer, Berlin 1889, 313-342 

\bibitem{dirichlet39}{\sc P.G.L. Dirichlet}, Recherches sur diverses applications de l'analyse infinit\'esimale a la th\'eorie des nombres, {\it J. reine angew. Math.} {\bf 19} (1839), 324-369; {\bf 21} (1840), 1-12, 134-155; Werke I, Reimer, Berlin 1889, 411-496 

\bibitem{diri}{\sc P.G.L. Dirichlet}, {\it Vorlesungen \"uber Zahlentheorie}, 3rd ed. with supplements by {\sc R. Dedekind}, Vieweg 1879

\bibitem{elstrodt}{\sc J. Elstrodt}, The Life and Work of Gustav Lejeune Dirichlet (1805-1859), in: {\it Analytic Number Theory: A Tribute to Gauss and Dirichlet}, W. Duke et al. (eds.), {\it Clay Math. Proc.} {\bf 7} (2007), 1-37

\bibitem{epstein}{\sc P. Epstein}, Zur Theorie allgemeiner Zetafunctionen, {\it Math. Ann.} {\bf 56} (1903), 615-644

\bibitem{eulerbern}{\sc L. Euler}, Demonstration de la somme de cette suite $1+1/4+1/9+1/16+\ldots$, {\it Journal Lit. d'Allemagne, de Suisse et du Nord} {\bf 2} (1743), 115-127; {\it Opera Omnia} {\bf I.14}, Teubner 1924, 138-155

\bibitem{eulerprod}{\sc L. Euler}, Variae observationes circa series infinitas, {\it Comment. Acad. Sci. Petropol} {\bf 12} (1737/1744), 53-96; {\it Opera Omnia} {\bf I.14}, Teubner 1924, 407-462

\bibitem{eulerfunc}{\sc L. Euler}, Remarques sur un beau rapport entre les s\'eries des puissances tant directes qui r\'eciproques, {\it M\'emoires de l'Acad\'emie des Sciences de Berlin} {\bf 17} (1749/1768), 83-106; {\it Opera Omnia} {\bf I.15}, Teubner 1924, 70-90

\bibitem{fine}{\sc N.J. Fine}, Note on the Hurwitz Zeta-Function, {\it Proc. A.M.S.} {\bf 2} (1951), 361-364

\bibitem{stamm}{\sc G. Frei, U. Stammbach}, {\it Die Mathematiker an den Z\"urcher Hochschulen}, Birkh\"auser, Basel 1994 

\bibitem{frobenius}{\sc F.G. Frobenius}, \"Uber Gruppencharaktere, {\it Berl. Ber.} (1896), 985-1021

\bibitem{garunklauri}{\sc R. Garunk\v stis, A. Laurin\v cikas}, {\it The Lerch Zeta-function}, Kluwer, Dordrecht 2002

\bibitem{gauss}{\sc C.F. Gauss}, {\it Werke, I}, G\"ottingen, 1870

\bibitem{gelb}{\sc St.S. Gelbart, St.D. Miller}, Riemann's Zeta Function and Beyond, {\it Bull. A.M.S.} {\bf 41} (2004), 59-112

\bibitem{grommer}{\sc J. Grommer}, Ganze transzendente Funktionen mit lauter reellen Nullstellen, {\it J. reine angew. Math.} {\bf 144} (1914), 114-165 

\bibitem{haas}{\sc R. Haas}, Goldbach, Hurwitz, and the Infinitude of Primes: Weaving a Proof across the Centuries, {\it Math. Intelligencer} {\bf 36} (2013), 54-60 

\bibitem{hada92}{\sc J. Hadamard}, Essai sur l'\'etude des fonctions donn\'ees par leur d\'eveloppement de Taylor, {\it Journ. de Math.} {\bf VIII} (1892), 101-186

\bibitem{hada93}{\sc J. Hadamard}, \'Etude sur le propri\'et\'es des fonctions enti\`eres et en particulier d'une fonction consid\'er\'ee par Riemann, {\it J. math. pures appl.} {\bf 9} (1893), 171-215 

\bibitem{hada962}{\sc J. Hadamard}, Sur les z\'eros de la fonction $\zeta(s)$ de Riemann, {\it C. R.} {\bf 122} (1896), 1470-1473 

\bibitem{hada96}{\sc J. Hadamard}, Sur la distribution des z\'eros de la fonction $\zeta(s)$ et ses cons\'equences arithm\'etiques, {\it Bull. Soc. Math. France} {\bf 24} (1896), 199-220 

\bibitem{hamburger}{\sc H. Hamburger}, \"Uber einige Beziehungen, die mit der Funktionalgleichung der Riemannschen $\zeta$-Funktion \"aquivalent sind, {\it Math. Ann.} {\bf 85} (1922), 129-140 

\bibitem{hankel}{\sc H. Hankel}, Die Euler'schen Integrale bei unbeschr\"ankter Variabilit\"at der Arguments, Habilitations-Diss., Leipzig 1863

\bibitem{hardy}{\sc G.H. Hardy}, Sur les z\'eros de la fonction $\zeta(s)$ de Riemann,
{\it Comptes Rendus Acad. Sci. Paris} {\bf 158} (1914), 1012-1014

\bibitem{hasse}{\sc H. Hasse}, Ein Summierungsverfahren f\"ur die Riemannsche $\zeta$-Reihe, {\it Math. Z.} {\bf 32} (1930), 458-464

\bibitem{hawk}{\sc St. Hawking}, Zeta function regularization of path integrals in curved spacetime, {\it Comm. Math. Phys.} {\bf 55} (1977), 133-148

\bibitem{hecke0}{\sc E. Hecke}, \"Uber die Zetafunktionen beliebiger algebraischer Zahlk\"orper, {\it Nachr. Ges. Wiss. G\"ottingen} (1917), 77-89

\bibitem{hecke}{\sc E. Hecke}, \"Uber die Bestimmung Dirichletscher Reihen durch ihre Funktionalgleichung, {\it Math. Ann.} {\bf 112} (1936), 664-699 

\bibitem{hejhal}{\sc D. Hejhal}, On a result by G. P\'olya concnerning the Riemann $\xi$-function, {\it J. Analyse Math.} {\bf 55} (1990), 59-95

\bibitem{hurwzeta}{\sc A. Hurwitz}, Einige Eigenschaften der Dirichletschen Funktionen $F(s)=\sum \left({D\over n}\right){1\over n^s}$, die bei der Bestimmung der Klassenzahlen bin\"arer quadratischer Formen auftreten, {\it Zeitschrift f. Math. u. Physik} {\bf 27} (1882), 86-101

\bibitem{hurwitz90}{\sc A. Hurwitz}, Ueber die Nullstellen der Bessel'schen Function, {\it Math. Annalen} {\bf 33} (1889), 246-266 

\bibitem{hurwreport1}{\sc A. Hurwitz}, Review to \cite{cahen2}, {\it Jahresberichte} JFM 25.0702.01

\bibitem{hurwitz14}{\sc A. Hurwitz}, Sur les points critiques des fonctions inverses dse fonctions enti\`eres, {\it Comptes Rendus} {\bf 158} (1914), 1007-1008

\bibitem{hurwitzdiaries}{\sc A. Hurwitz}, {\it Die Mathematischen Tageb\"ucher und der \"ubrige handschriftliche Nachlass von Adolf Hurwitz}, Handschriften und Autographen der ETH Z\"urich, {\sf http://www.e-manuscripta.ch/} 

\bibitem{hurwitzbuch}{\sc A. Hurwitz}, {\it Lectures on Number Theory}, translated from unpublished German course notes by {\sc N. Kritikos}, Springer 1985

\bibitem{iwaniec}{\sc H. Iwaniec}, {\it Lectures on the Riemann zeta function}, AMS, Providence 2014 

\bibitem{jv}{\sc J.L.W.V. Jensen}, Unders{\o}gelser over Ligningernes Theori, {\it Beretning om den anden skandinaviske Matematikerkongres i Kj{\o}benhavn} 1911, 51-65 

\bibitem{jonquiere}{\sc A. Jonqui\`ere}, Note sur la s\'erie $\sum_{n=1}^\infty {x^n\over n^s}$, {\it Bull. de la S.M.F.} {\bf 17} (1889), 142-152

\bibitem{kinkelin}{\sc H. Kinkelin}, Allgemeine Theorie der harmonischen Reihen, mit Anwendungen auf die Zahlentheorie, {\it Programm der Gewerbeschule Basel} (1861/62), 1-32

\bibitem{koch}{\sc H. Koch}, Oskar Xaver Schl\"omilch -- ein F\"orderer des mathematischen Unterrichts f\"ur Techniker und Ingenieure, {\it NTM, Schriftenr. Gesch. Naturwiss. Tech. Med.} {\bf 27} (1990), 1-10

\bibitem{kronecker}{\sc L. Kronecker}, Notiz \"uber Potenzreihen, {\it Monatsber. Preuss. Akad. Wiss. Berlin} (1878), 53-58; Werke {\bf V}, Teubner, Leipzig 1930, 197-201

\bibitem{landau03}{\sc E. Landau}, Neuer Beweis des Primzahlsatzes und Beweis des Primidealsatzes, {\it Math. Annalen} {\bf 56} (1903), 645-670

\bibitem{landau}{\sc E. Landau}, Euler und die Funktionalgleichung der Riemannschen Zetafunktion, {\it Bibl. Mat.} {\bf 7} (1906), 69-79

\bibitem{landau09}{\sc E. Landau}, {\it Handbuch der Lehre von der Verteilung der Primzahlen} in zwei B\"anden, Teubner, Leipzig 1909

\bibitem{landau27}{\sc E. Landau}, {\it Vorlesungen \"uber Zahlentheorie}, Zweiter Band, Hirzel, Leipzig 1927

\bibitem{lang0}{\sc S. Lang}, {\it Introduction to Modular Forms}, Springer 1976

\bibitem{lang}{\sc S. Lang}, {\it Complex Analysis}, Springer, 1999, 4th ed.

\bibitem{laugw}{\sc D. Laugwitz}, {\it Bernhard Riemann 1826-1866}, Birkh\"auser Boston, 2008, engl. translation 

\bibitem{legend}{\sc A.M. Legendre}, {\it Essai sur la th\'eorie des nombres}, Courcier Paris, 2nd ed. 1808

\bibitem{lerch87}{\sc M. Lerch}, Sur la fonction ${\mathcal K}(w,x,s)=\sum_{k=0}^\infty {\exp(2k\pi ix)\over (w+k)^s}$, {\it Acta Math.} {\bf 11} (1887/1888), 19-24

\bibitem{lerch92}{\sc M. Lerch}, Z\'akladove theorie Malmst\`enovych \v rad, {\it Rozpravy \v Cesk\`e Akad.} 2.Kl. {\bf 1} (1892), 1-70

\bibitem{lindeloef}{\sc E. Lindel\"of}, Robert Hjalmar Mellin, {\it Acta Math.} {\bf 61} (1933), i-vii

\bibitem{lipschitz}{\sc R. Lipschitz}, Untersuchung einer aus vier Elementen gebildeten Reihe, {\it J. reine angew. Math.} {\bf 54} (1857), 313-328

\bibitem{lip89}{\sc R. Lipschitz}, Untersuchungen der Eigenschaften einer Gattung von unendlichen Reihen, {\it J. reine angew. Math.} {\bf 105} (1889), 127-156

\bibitem{malmsten}{\sc C.J. Malmst\`en}, De integralibus quibusdam definitis, seriebusque infinitis, {\it J. reine angew. Mathe.} {\bf 38} (1849), 1-39\footnote{An English translation due to A. Aycock is available as {\sf arXiv:1309.3824v1} in the {\sf arXiv}.}

\bibitem{mango95}{\sc H. von Mangoldt}, Zu Riemann's Abhandlung ''Ueber die Anzahl der Primzahlen unter einer gegebenen Gr\"o\ss e'', {\it J. reine angew. Math.} {\bf 114} (1895), 255-305

\bibitem{mellin1}{\sc Hj. Mellin}, \"Uber eine Verallgemeinerung der Riemannschen Function $\zeta(s)$, {\it Acta soc. scientiarium Fennicae} {\bf 24} (1899), no. 10, 50 p

\bibitem{mellin2}{\sc Hj. Mellin}, Die Dirichlet'schen Reihen, die zahlentheoretischen Funktionen und die unendlichen Produkte von endlichem Geschlecht, {\it Acta Math.} {\bf 28} (1904), 37-64

\bibitem{mertens}{\sc F. Mertens}, Beweis, dass jede linerae Function mit ganzen complexen teilerfremden Coefficienten unendlich viele complexe Primzahlen darstellt, {\it Wien Ber.} {\bf 108} (1899), 517-556 

\bibitem{minkowski}{\sc H. Minkowski}, Peter Gustav Lejeune Dirichlet und seine Bedeutung f\"ur die heutige Mathematik, {\it Jber. DMV} {\bf 14} (1905), 149-163

\bibitem{zass}{\sc H. Minkowski}, {\it Briefe an David Hilbert}, {\sc L. R\"udenberg, H. Zassenhaus} (eds.), Springer 1973

\bibitem{murty}{\sc M.R. Murty, N. Thain}, Prime Numbers in Certain Arithmetic Progressions, {\it Funct. Approximatio. Comment. Math.} {\bf 35} (2006), 249-259

\bibitem{narkie}{\sc W. Narkiewicz}, {\it The Development of Prime Number Theory}, Springer 2000

\bibitem{oswald}{\sc N.M.R. Oswald, J. Steuding}, Complex continued fractions: early work of the brothers Adolf and Julius Hurwitz, {\it Arch. Hist. Exact Sci.} {\bf 68} (2014), 499-528

\bibitem{perron}{\sc O. Perron}, Zur Theorie der Dirichletschen Reihen, {\it J. reine angew. Math.} {\bf 134} (1908), 95-143

\bibitem{poisson}{\sc S. Poisson}, {\it Theorie math\'ematiques de la chaleur}, Bachelier, Paris 1835 

\bibitem{polya13}{\sc G. P\'olya}, \"Uber Ann\"aherung durch Polynome mit lauter reellen Wurzeln, {\it Rend. Circ. Mat. Palermo} {\bf 36} (1913), 279-295

\bibitem{polya13a}{\sc G. P\'olya}, Sur une question concernant les fonctions enti\`eres, {\it Comptes Rendus} {\bf 158} (1914), 330-333

\bibitem{polya18}{\sc G. P\'olya}, Ueber die Nullstellen gewisser ganzer Funktionen, {\it Math. Zeit.} {\bf 2} (1918), 352-383

\bibitem{polya26a}{\sc G. P\'olya}, Bemerkung \"uber die Integraldarstellung der Riemannschen $\xi$-Funktion, {\it Acta Math.} {\bf 48} (1926), 305-317

\bibitem{polya26b}{\sc G. P\'olya}, On the zeros of certain trigonometric integrals, {\it J. London Math. Soc.} {\bf 1} (1926), 98-99

\bibitem{polya27}{\sc G. P\'olya}, \"Uber die algebraisch-funktionentheoretischen Untersuchungen von J. L. W. Jensen, {\it Meddelelser Kobenhavn} {\bf 7} (1927), 3-33

\bibitem{polya27b}{\sc G. P\'olya}, \"Uber trigonometrische Integrale mit nur reellen Nullstellen, {\it } {\bf 158} (1927), 6-18

\bibitem{diaries}{\sc G. P\'olya}, \"Uber den handschriftlichen Nachlass des Mathematikers Adolf Hurwitz, Handschriften und Autographen der ETH-Bibliothek, 

\bibitem{odlyzkopolya}{\sc G. P\'olya}, Letter exchange between P\'olya and Odlyzko, (1982), available at {\sf www.dtc.umn.edu/$\sim$odlyzko/polya/}

\bibitem{posz}{\sc G. P\'olya, G. Szeg\"o}, {\it Aufgaben und Lehrs\"atze aus der Analysis} I, II, Springer 1925

\bibitem{porubsky}{\sc \v St. Porubsk\'y}, Maty\'a\v s Lerch (1860-1922), {\it \v Siauliai Math. Semin.} {\bf 8 (16)} (2013), 197-222

\bibitem{rie}{\sc B. Riemann}, \"Uber die Anzahl der Primzahlen unterhalb einer gegebenen Gr\"osse, {\it Monatsber. Preuss. Akad. Wiss. Berlin} (1859), 671-680

\bibitem{rowe0}{\sc D.E. Rowe}, Felix Klein, Adolf Hurwitz, and the ''Jewish Question'' in German Academia, {\it Math. Intelligencer} {\bf 29} (2007), 18-30

\bibitem{rowe}{\sc D.E. Rowe, R. Schulmann}, General Relativity in the Context of Weimar Culture, {\it Max-Planck-Institut Wissenschaftsgeschichte} (2014), preprint 456 

\bibitem{ida}{\sc I. Samuel-Hurwitz}, Erinnerungen an die Familie Hurwitz, mit Biographie ihres Gatten Adolph Hurwitz, Prof. f. h\"ohere Mathematik an der ETH, HS 583a: 2, Archive ETH Zurich University Archives

\bibitem{schaertlin}{\sc G. Schaertlin}, Hermann Kinkelin, {\it Mitteilungen schweizerischer versicherungsmathematiker} {\bf 28} (1933), 1-17

\bibitem{scharlau}{\sc W. Scharlau}, The mathematical correspondence of Rudolf Lipschitz, {\it Hist. Math.} {\bf 13} (1986), 165-167

\bibitem{scheibner}{\sc H. Scheibner}, Ueber die Anzahl der Primzahlen unter einer beliebigen Grenze, {\it Zeitschr. f. Math.} {\bf 5} (1860) 

\bibitem{schloemilch1}{\sc O. Schl\"omilch}, Uebungsaufgaben f\"ur Sch\"uler ueber den Lehrsatz von dem Herrn Prof. Dr. Schl\"omilch, {\it Arch. Math. Phys.} {\bf 12} (1849), 415

\bibitem{schloemilch2}{\sc O. Schl\"omilch}, \"Uber eine Eigenschaft gewisser Reihen, {\it Zeitschr. Math. Phys.} {\bf 3} (1858), 130-132

\bibitem{schloemilch278}{\sc O. Schl\"omilch}, Ueber die Summen von Potenzen der reciproken nat\"urlichen Zahlen, {\it Zeitschr. Math. Phys.} {\bf 23} (1878), 135-137

\bibitem{schnee}{\sc W. Schnee}, Die Funktionsgleichung der Zetafunktion und der Dirichletschen Reihen mit periodischen Koeffizienten, {\it Math. Z.} {\bf 31} (1929), 378-390

\bibitem{wschwarz}{\sc W. Schwarz}, Some Remarks on the History of the Prime Number Theorem from 1896 to 1960, in: {\it Development of Mathematics 1900-1950}, edited by {\sc J.-P. Pier}, Birkh\"auser 1994, 565-616

\bibitem{selberg}{\sc A. Selberg}, Old and new conjectures and results about a class of Dirichlet series, in: `{\it Proceedings of the Amalfi Conference on Analytic Number Theory}', Maiori 1989, E. Bombieri et al. (eds.), Universit\`a di Salerno 1992, 367-385

\bibitem{siegel}{\sc C.L. Siegel}, \"Uber Riemanns Nachla\ss{} zur analytischen Zahlentheorie, {\it Quell. Stud. Gesch. Mat. Astr. Physik} {\bf 2} (1932), 45-80

\bibitem{siegel44}{\sc C.L. Siegel}, The average measure of quadratic forms with given determinant and signature, {\it Ann. of Math.} {\bf 45} (1944), 667-685

\bibitem{siegelh}{\sc C.L. Siegel}, Zur Geschichte des Frankfurter Mathematischen Seminars, in: {\it Gesammelte Abhandlungen, Vol. III}, Springer 1966

\bibitem{tate}{\sc J. Tate}, Fourier analysis in number fields, and Hecke's zeta-functions, in: {\it Algebraic Number Theory}, Proc. Instructional Conf. Brighton 1965, Thompson, Washingthon 1967, 305-347

\bibitem{titch}{\sc E.C. Titchmarsh}, On Epstein's zeta-function, {\it Proc. London Math. Soc.} {\bf 36} (1934), 485-500

\bibitem{poussin96}{\sc C.J. de la Vall\'ee Poussin}, Recherches analytique sur la th\'eorie des nombres premiers, I-III, {\it Ann. Soc. Sci. Bruxelles} {\bf 20} (1896), 183-256, 281-362, 363-397

\bibitem{vapo}{\sc C.J. de la Vall\'ee Poussin}, Recherches analytique sur la th\'eorie des nombres premiers, IV,V, {\it Ann. Soc. Sci. Bruxelles} {\bf 21} (1897), 251-342, 343-368

\bibitem{weil}{\sc A. Weil}, Prehistory of the Zeta-Function, in: {\it Number theory, trace formulas and discrete groups} (Oslo, 1987) (ed. E.A. Aubert et al.), Academic Press, Boston, MA, 1989, 1-9 

\bibitem{weil2}{\sc A. Weil}, On Eisenstein's Copy of the Disquisitiones, in: {\it Algebraic number theory -- in honor of K. Iwasawa}, Proc. Workshop Iwasawa Theory Spec. Values $L$-Funct., Berkeley/CA (USA) 1987, {\it Adv. Stud. Pure Math.} {\bf 17} (1989), 463-469 

\end{thebibliography}
\end{document}